\documentclass[aos,MSNbibl,seceqn,nameyear,rotating,dvips]{arximspdf}
\usepackage{dcolumn}
\usepackage{graphicx}

\doi{10.1214/12-AOS1062} 
\volume{41}
\issue{1}
\pubyear{2013}
\firstpage{1}
\lastpage{40}

\makeatletter

\newcolumntype{d}[1]{D{.}{.}{#1}}

\renewcommand{\widehat}{\hat}

\newcommand{\uo}{\mathrm{uo}}

\newcommand{\st}{\mathrm{st}}

\newcommand{\lst}{l\mbox{-}\mathrm{st}}

\renewcommand{\overline}{\bar}
\newcommand{\rrVert}{\Vert}
\newcommand{\rrvert}{\vert}
\newcommand{\llVert}{\Vert}
\newcommand{\llvert}{\vert}

\newtheorem{theorem}{Theorem}[section]
\newproclaim{algorithm}{Algorithm}[section]
\newproclaim{definition}{Definition}[section]
\newproclaim{example}{Example}[section]
\newtheorem{lemma}{Lemma}[section]
\newtheorem{proposition}{Proposition}[section]
\newproclaim{remark}{Remark}[section]

\makeatother

\begin{document}
\begin{frontmatter}

\title{The linear stochastic order and directed inference for multivariate
ordered distributions}
\runtitle{Inference for the linear stochastic order}

\begin{aug}
\author[A]{\fnms{Ori} \snm{Davidov}\corref{}\thanksref{t1}\ead[label=e1]{davidov@stat.haifa.ac.il}}
\and
\author[B]{\fnms{Shyamal} \snm{Peddada}\thanksref{t2}\ead[label=e2]{peddada@niehs.nih.gov}}
\runauthor{O. Davidov and S. Peddada}
\affiliation{University of Haifa and National Institute of
Environmental Health~Sciences}
\address[A]{Department of Statistics\\
University of Haifa\\
Mount Carmel, Haifa 31905\\
Israel\\
\printead{e1}} 
\address[B]{Biostatistics Branch\\
National Institute\\
\quad of Environmental Health Sciences\\
111 T.W. Alexander Drive\\
Research Triangle Park, North Carolina 27709\\
USA\\
\printead{e2}}
\end{aug}

\thankstext{t1}{Supported in part by the Israeli Science
Foundation Grant No 875/09.}

\thankstext{t2}{Supported in part by the Intramural
Research Program of the NIH, National Institute of Environmental Health
Sciences (Z01 ES101744-04).}

\received{\smonth{8} \syear{2011}}
\revised{\smonth{7} \syear{2012}}

%
\begin{abstract}
Researchers are often interested in drawing inferences regarding the
order between two experimental groups on the basis of multivariate
response data. Since standard multivariate methods are designed for
two-sided alternatives, they may not be ideal for testing for order
between two groups. In this article we introduce the notion of the
linear stochastic order and investigate its properties. Statistical
theory and methodology are developed to both estimate the direction
which best separates two arbitrary ordered distributions and to test
for order between the two groups. The new methodology generalizes Roy's
classical largest root test to the nonparametric setting and is
applicable to random vectors with discrete and/or continuous
components. The proposed methodology is illustrated using data obtained
from a 90-day pre-chronic rodent cancer bioassay study conducted by the
National Toxicology Program (NTP).
\end{abstract}

%
\begin{keyword}[class=AMS]
\kwd{60E15}
\kwd{62E20}
\kwd{62G10}
\kwd{62G20}
\kwd{62H99}
\kwd{62P15}
\end{keyword}
\begin{keyword}
\kwd{Nonparametric tests}
\kwd{order-restricted statistical inference}
\kwd{stochastic order relations}
\end{keyword}

\end{frontmatter}

\section{Introduction}\label{sec1}

In a variety of applications researchers are interested in comparing
two treatment groups on the basis of several, potentially dependent
outcomes. For example, to evaluate if a chemical is a neuro-toxicant,
toxicologists compare a treated group of animals with an untreated
control group in terms of various correlated outcomes such as
tail-pinch response, click response and gait score, etc.;
cf. \citet{Mos00}. The statistical problem of interest is to
compare the multivariate distributions of the outcomes in the control
and treatment groups. Moreover, the outcome distributions are expected
to be ordered in some sense. The theory of stochastic order relations
[\citet{ShaSha07}] provides the theoretical foundation for such
comparisons.

To fix ideas let $\mathbf{X}$ and $\mathbf{Y}$ be $p$-dimensional
random variables (RVs); $\mathbf{X}$ is said to be smaller than
$\mathbf{Y}$ in the multivariate stochastic order, denoted
$\mathbf{X}%
\preceq_{\st }\mathbf{Y}$, provided $\mathbb{P} ( \mathbf{X}\in U )
\leq\mathbb{P} ( \mathbf{Y}\in U ) $ for all upper sets $U\in
\mathbb{R}^{p}$ [\citet{ShaSha07}]. If for some upper set
the above inequality is sharp, we say that $\mathbf{X}$ is strictly
smaller than $\mathbf{Y}$ (in the multivariate stochastic order) which
we denote by $\mathbf{X}\prec_{\st }\mathbf{Y}$. Recall that a set
$U\in\mathbb{R}^{p}$ is called an upper set if $\mathbf{u}\in
U$ implies that $\mathbf{v}\in U$ whenever $\mathbf{u}\leq$ $\mathbf{v}$,
that is, if $u_{i}\leq v_{i}$, $i=1,\ldots,p$. Note that comparing
$\mathbf{X}$ and $\mathbf{Y}$ with respect to the multivariate stochastic
order requires comparing their distributions over all upper sets in $
\mathbb{R}^{p}$. This turns out to be a very high-dimensional problem.
For example, if $\mathbf{X}$ and $\mathbf{Y}$ are multivariate binary RVs,
then $\mathbf{X}\preceq_{\st }\mathbf{Y}$ provided $\sum_{\mathbf{t}\in
U}p_{\mathbf{X}} ( \mathbf{t} ) \leq\sum_{\mathbf{t}\in U}p_{\mathbf{Y}} (
\mathbf{t} ) $ where $p_{\mathbf
{X}%
} ( \mathbf{t} ) $ and $p_{\mathbf{Y}} ( \mathbf{t} ) $ are the
corresponding probability mass functions. Here $U\in\mathcal{U}_{p}$
where $\mathcal{U}_{p}$ is the family of upper sets defined on the
support of a $p$-dimensional multivariate binary RV. It turns out that
the cardinality of $\mathcal{U}_{p}$, denoted by $N_{p}$, grows
super-exponentially with $p$. In fact $N_{1}=1$, $N_{2}=4$, $N_{3}=18$,
$N_{4}=166$, $N_{5}=7579$ and $N_{6}=7\mbox{,}828\mbox{,}352$. The
values of $N_{7}$ and $N_{8}$ are also known, but $N_{9}$ is not.
However, good approximations for $N_{p}$ are available for all $p$; cf.
\citet{DavPed11}. Obviously the number of upper sets for general
multivariate RVs is much larger. Since in many applications $p$ is
large, it would seem that the analysis of high-dimensional
stochastically ordered data is practically hopeless. As a consequence,
the methodology for analyzing multivariate ordered data is
underdeveloped. It is worth mentioning that \citet{SamWhi89} as
well as \citet{LucWri91} studied stochastically ordered bivariate
multinomial distributions. They noted the difficulty of extending their
methodology to high-dimensional data due to the large number of
constraints that need to be imposed. Recently \citet{DavPed11}
proposed a framework for testing for order among $K$, $p$-dimensional,
ordered multivariate binary distributions.

In this paper we address the dimensionality problem by considering an
easy to
understand stochastic order which we refer to as the linear stochastic order.
%
\begin{definition}
\label{Def-lst}The RV $\mathbf{X}$ is said to be smaller than the RV
$\mathbf{Y}$ in the (multivariate) linear stochastic order, denoted
$\mathbf{X}\preceq_{\lst }\mathbf{Y}$, if for all $\mathbf
{s}\in$ $
\mathbb{R}_{+}^{p}= \{ \mathbf{s}\dvtx \mathbf{s}\geq0 \}, $
%
\begin{equation}\label{l-st}
\mathbf{s}^{T}\mathbf{X}\preceq_{\st }\mathbf{s}^{T}
\mathbf{Y},
\end{equation}
where $\preceq_{\st }$ in (\ref{l-st}) denotes the usual (univariate)
stochastic order.
\end{definition}

Note that it is enough to limit (\ref{l-st}) to all nonnegative real
vectors satisfying $\llVert\mathbf{s}\rrVert=1$, and
accordingly we denote by $\mathcal{S}_{+}^{p-1}$\vadjust{\goodbreak} the positive part of
the unit sphere in $ \mathbb{R}^{p}$. We call each
$\mathbf{s}\in\mathcal{S}_{+}^{p-1}$ a ``direction.'' In other words
the RVs $\mathbf{X}$ and $\mathbf{Y}$ are ordered by the linear
stochastic order if every nonnegative linear combination of their
components is ordered by the usual (univariate) stochastic order. Thus
instead of considering all upper sets in $ \mathbb{R}^{p}$ we need for
each $\mathbf{s}\in\mathcal{S}_{+}^{p-1}$ to consider only upper
sets in $\mathbb{R}$. This is a substantial reduction in
dimensionality. In fact we will show that only one value of
$\mathbf{s}$ need be considered. Note that
the linear stochastic order, like the multivariate stochastic order, is
a generalization of the usual univariate stochastic order to
multivariate data. Both of these orders indicate, in different ways,
that one random vector is more likely than another to take on large
values. In this paper we develop the statistical theory and methodology
for estimation and testing for linearly ordered multivariate
distributions. For completeness we note that weaker notions of the
linear stochastic order are discussed by \citet{HuHomMeh} and applied
to various optimization problems in queuing and finance.

Comparing linear combinations has a long history in statistics. For example,
in Phase I clinical trials it is common to compare dose groups using an
overall measure of toxicity. Typically, this quantity is an ad hoc weighted
average of individual toxicities where the weights are often known as
``severity weights;'' cf. \citet{BekTha04} and \citet{IvaMur09}.
This strategy of dimension reduction is not new in the statistical literature
and has been used in classical multivariate analysis when comparing two or
more multivariate normal populations. For example, using the
union-intersection principle, the comparison of multivariate normal
populations can be reduced to the comparison of all possible linear
combinations of their mean vectors. This approach is the basis of Roy's
classical largest root test [\citet{Roy53}, \citet{JohWic98}]. Our
proposed test may be viewed as nonparametric generalization of the classical
normal theory method described above with the exception that we limit
consideration only to nonnegative linear combinations (rather than all
possible linear combinations) since our main focus is to make
comparisons in
terms of stochastic order. We emphasize that the linear stochastic
order will
allow us to address the much broader problem of directional ordering for
multivariate ordered data, that is, to find the direction which best separates
two ordered distributions. Based on our survey of the literature, we
are not
aware of any methodology that addresses the problems investigated here.

This paper is organized in the following way. In Section~\ref{sec2} some probabilistic
properties of the linear stochastic order are explored, and its relationships
with other multivariate stochastic orders are clarified. In Section~\ref{sec3} we
provide the background and motivation for directional inference under the
linear stochastic order and develop estimation and testing procedure for
independent as well as paired samples. In particular the estimator of
the best
separating direction is presented and\vadjust{\goodbreak} its large sampling properties derived.
We note that the problem of estimating the best separating direction is a
nonsmooth optimization problem. The limiting distribution of the best
separating direction is derived in a variety of settings. Tests for the linear
stochastic order based on the best separating direction are also developed.
One advantage of our approach is that it avoids the estimation of multivariate
distributions subject to order restrictions. Simulation results,
presented in
Section~\ref{sec4}, reveal that for large sample sizes the proposed estimator has
negligible bias and mean squared error (MSE). The bias and MSE seem to depend
on the true value of the best separating direction, the dependence structure
and the dimension of the problem. Furthermore, the proposed test honors the
nominal type I error rate and has sufficient power. In Section~\ref{sec5} we illustrate
the methodology using data obtained from the National Toxicology Program
(NTP). Concluding remarks and some open research problems are provided in
Section~\ref{sec6}. For convenience all proofs are provided in the
\hyperref[app]{Appendix} where
additional concepts are defined when needed.

\section{Some properties of the linear stochastic order}\label{sec2}

We start by clarifying the relationship between the linear stochastic order
and the multivariate stochastic order. First note that $\mathbf{X}%
\preceq_{\lst }\mathbf{Y}$ if and only if $\mathbb{P}(\mathbf{s}%
^{T}\mathbf{X}\geq t)\leq\mathbb{P}(\mathbf{s}^{T}\mathbf
{Y}\geq t)$ for all $(t,\mathbf{s})\in\mathbb{R} \times\mathbb
{R}_{+}^{p}$ which is equivalent to $\mathbb{P}(\mathbf{X}\in H)\leq
\mathbb{P}(\mathbf{Y}\in H)$ for all $H\in\mathcal{H}$ where
$\mathcal{H}$ is the collection of all upper half-planes, that is, sets
which are both half planes and upper sets. Thus
$\mathbf{X}\preceq_{\st }\mathbf{Y} $ $\Rightarrow$
$\mathbf{X}\preceq_{\lst }\mathbf{Y}$. The converse does not
hold in general.
%
\begin{example}
Let $\mathbf{X}$ and $\mathbf{Y}$ be bivariate RVs such that
$\mathbb{P}(\mathbf{X}= ( 1,1 ) )=\mathbb{P}(\mathbf{X}
= ( 0,1 ) )=\mathbb{P}(\mathbf{X}= ( 1,0 ) )=1/3$
and $\mathbb{P}(\mathbf{Y}= ( 3/4,3/4 ) )=\mathbb{P}(%
\mathbf{Y}= ( 1,2 ) )=\mathbb{P}(\mathbf{Y}= (
2,1 ) )=1/3$. It is easy to show that $\mathbf{X}$ is smaller than
$\mathbf{Y}$ in the linear stochastic order but not in the multivariate
stochastic order.
\end{example}

The following theorem provides some closure results for the linear
stochastic order.
%
\begin{theorem}
\label{Thm-Closure} \textup{(i)} If $\mathbf{X}\preceq_{\lst
}\mathbf{Y}$, then $g(\mathbf{X})\preceq_{\lst }g(\mathbf{Y})$ for any
affine increasing function; \textup{(ii)} if $\mathbf
{X}%
\preceq_{\lst }\mathbf{Y}$, then $\mathbf{X}_{I}\preceq_{\lst }%
\mathbf{Y}_{I}$ for each subset $I\in\{ 1,\ldots,p \};$
\textup{(iii)}~if $\mathbf{X}|\mathbf{Z}=\mathbf{z}\preceq_{\lst
}\mathbf{Y}|\mathbf{Z}=\mathbf {z}$ for all $\mathbf{z}$ in the support
of $\mathbf{Z}$, then $\mathbf{X}\preceq_{\lst }\mathbf{Y};$
\textup{(iv)} if $\mathbf{X}_{1},\ldots,\mathbf{X}_{n}$ are
independent RVs with dimensions $p_{i}$ and similarly for $\mathbf
{Y}_{1},\ldots,\mathbf{Y}_{n}$ and if in addition $\mathbf{X}_{i}%
\preceq_{\lst }\mathbf{Y}_{i}$, then
$(\mathbf{X}_{1},\ldots,\mathbf{X}_{n})\preceq_{\lst
}(\mathbf{Y}_{1},\ldots,\mathbf {Y}_{n})$; \textup{(v)} finally, if
$\mathbf{X}_{n}\rightarrow \mathbf{X}$ and $\mathbf{Y}_{n}$
$\rightarrow\mathbf{Y}$ where convergence can be in distribution, in
probability or almost surely and if $\mathbf{X}_{n}\preceq_{\lst
}\mathbf{Y}_{n}$ for all $n$, then $\mathbf{X}\preceq_{\lst
}\mathbf{Y}$.
\end{theorem}

Theorem~\ref{Thm-Closure} shows that the linear stochastic order is closed
under increasing linear transformations, marginalization, mixtures,
conjugations and convergence. In particular parts\vadjust{\goodbreak} (ii) and
(iii) of Theorem~\ref{Thm-Closure} imply that if
$\mathbf{X}\preceq_{\lst }\mathbf{Y}$, then $X_{i}\preceq_{\st }Y_{i}$ and
$X_{i}+X_{j}\preceq_{\st }Y_{i}+Y_{j}$ for all $i$ and $j$; that is, all marginals
are ordered as are all convolutions. Although the multivariate stochastic
order is in general stronger than the linear stochastic order, there are
situation in which both orders coincide.
%
\begin{theorem}
\label{Thm-ERVs}Let $\mathbf{X}$ and $\mathbf{Y}$ be continuous
elliptically distributed RVs supported on $\mathbb{R}^{p}$ with the
same generator. Then $\mathbf{X}\preceq_{\lst }%
\mathbf{Y}$ if and only if $\mathbf{X}\preceq_{\st }\mathbf{Y}$.
\end{theorem}

Note that the elliptical family of distributions is large and includes the
multivariate normal, multivariate $t$ and the exponential power family; see
\citet{FanKotNg89}. Thus Theorem~\ref{Thm-ERVs} shows that the multivariate
stochastic order coincides with the linear stochastic order in the normal
family. Incidentally, in the proof of Theorem~\ref{Thm-ERVs} we
generalize the
results of \citet{DinZha04} on multivariate stochastic ordering of
elliptical RVs. Another interesting example is the following:
%
\begin{theorem}
\label{Thm-MVBs}Let $\mathbf{X}$ and $\mathbf{Y}$ be multivariate
binary RVs. Then $\mathbf{X}\preceq_{\lst }\mathbf{Y}$ is
equivalent to
$\mathbf{X}\preceq_{\st }\mathbf{Y}$ if and only if $p\leq3$.
\end{theorem}
%
\begin{remark}
In the proof of Theorem~\ref{Thm-ERVs} distributional properties of the
elliptical family play a major role. In contrast, Theorem \ref
{Thm-MVBs} is a
consequence of the geometry of the upper sets of multivariate binary
RVs which
turn out to be upper half planes if and only if $p\leq3$.
\end{remark}

We now explore the role of the dependence structure.
%
\begin{theorem}
\label{Thm-Copula}Let $\mathbf{X}$ and $\mathbf{Y}$ have the same
copula. Then $\mathbf{X}\preceq_{\lst }\mathbf{Y}$ if and only if
$\mathbf{X}\preceq_{\st }\mathbf{Y}$.
\end{theorem}

Theorem~\ref{Thm-Copula} establishes that if two RVs have the same dependence
structure as quantified by their copula function [cf. \citet{Joe97}], then the
linear and multivariate stochastic orders coincide. Such situations
arise when
the correlation structure among outcomes is not expected to vary with dose.

The orthant orders are also of interest in statistical applications. We
say that $\mathbf{X}$ is smaller than $\mathbf{Y}$ in the upper orthant
order, denoted $\mathbf{X}\preceq_{\uo}\mathbf{Y}$, if $\mathbb
{P}(\mathbf{X}\in O)\leq\mathbb{P}(\mathbf{Y}\in O)$ for all
$O\in\mathcal{O}$ where $\mathcal{O}$ is the collection of upper
orthants, that is, sets of the form $ \{ \mathbf{z}\dvtx \mathbf{z}\geq
\mathbf{x} \} $ for some fixed $\mathbf{x}\in \mathbb{R}^{p}$. The
lower orthant order is similarly defined; cf. \citet{ShaSha07} or
\citet{DavHer11}. It is obvious that the orthant orders are weaker
than the usual multivariate stochastic order, that is,
$\mathbf{X}\preceq_{\st }\mathbf{Y}$ $\Rightarrow$ $\mathbf{X} 
\preceq_{\uo}\mathbf{Y}$ and $\mathbf{X}\preceq_{\mathrm{lo}}\mathbf
{Y}$. In
general the linear stochastic order does not imply the upper (or lower)
orthant order, nor is the converse true. However, as stated below,
under some
conditions on the copula functions, the linear stochastic order implies the
upper (or lower) orthant order.\vadjust{\goodbreak}
%
\begin{theorem}
\label{Thm-ORT}If $\mathbf{X}\preceq_{\lst }\mathbf{Y}$ and
$C_{\mathbf{X}} ( \mathbf{u} ) \leq C_{\mathbf{Y}%
} ( \mathbf{u} ) $ for all $\mathbf{u}\in[0,1]^{p}$,
then $\mathbf{X}\preceq_{\mathrm{lo}}\mathbf{Y}$. Similarly if
$\mathbf{X}%
\preceq_{\lst }\mathbf{Y}$ and $\overline{C}_{\mathbf{X}} (
\mathbf{u} ) \leq\overline{C}_{\mathbf{Y}} (
\mathbf{u}%
) $ for all $\mathbf{u}\in[0,1]^{p}$, then $\mathbf
{X}%
\preceq_{\uo}\mathbf{Y}$.
\end{theorem}

Note that $C_{\mathbf{X}} ( \mathbf{u} ) $ and
$\overline{C}_{\mathbf{X}} ( \mathbf{u} ) $ above are the
copula and tail-copula functions for the RV $\mathbf{X}$ [cf. Joe (\citeyear{Joe97})]
and are defined in the \hyperref[app]{Appendix} and similarly for $C_{\mathbf
{Y}} (
\mathbf{u} ) $ and $\overline{C}_{\mathbf{Y}} (
\mathbf{u} )$. Further note that the relations $C_{\mathbf
{X}%
} ( \mathbf{u} ) \leq C_{\mathbf{Y}} ( \mathbf
{u}%
) $ and/or $\overline{C}_{\mathbf{X}} ( \mathbf
{u} )
\leq\overline{C}_{\mathbf{Y}} ( \mathbf{u} ) $ indicate
that the components of $\mathbf{Y}$ are more strongly dependent
than the
components of $\mathbf{X}$. This particular dependence ordering is known
as positive quadrant dependence. It can be further shown that strong
dependence and the linear stochastic order do not in general imply
stochastic ordering.

Additional properties of the linear stochastic order as they relate to
estimation and testing problems are given in Section
\ref{Sec-Dir}.

\section{Directional inference}\label{Sec-Dir}\label{sec3}

\subsection{Background and motivation}\label{sec3.1}

There exists a long history of well-deve\-loped theory for comparing two
or more multivariate normal (MVN) populations. Methods for assessing
whether there are any differences between the populations [which
differ? in which component(s)? and by how much?] have been addressed in
the literature using a variety of simultaneous confidence intervals and
multiple comparison methods; cf. \citet{JohWic98}. Of
particular interest to us is Roy's largest root test. To fix ideas
consider two multivariate normal random vectors $\mathbf{X}$ and
$\mathbf{Y}$ with means $\bolds{\mu}$ and $\bolds{\nu}$,
respectively, and a common variance matrix~$\bolds{\Sigma}$. Using
the union-intersection principle \citet{Roy53} expressed the problem of
testing $H_{0}\dvtx \bolds{\mu}=\bolds{\nu}$ versus
$H_{1}\dvtx \bolds{\mu}\neq\bolds{\nu}$ as a collection of
univariate testing problems, by showing that $H_{0}$ and $H_{1}$ are
equivalent to $\bigcap_{\mathbf{s}\in\mathbb{R}
^{p}}H_{0,\mathbf{s}}$ and $\bigcup_{\mathbf{s}\in
\mathbb{R}
^{p}}H_{1,\mathbf{s}}$ where $H_{0,\mathbf{s}}\dvtx \mathbf{s}^{T}%
\bolds{\mu}=\mathbf{s}^{T}\bolds{\nu}$ and
$H_{1,\mathbf{s}%
}\dvtx \mathbf{s}^{T}\bolds{\mu}\neq\mathbf{s}^{T}\bolds
{\nu}$.
Implicitly Roy's test identifies the linear combination $\mathbf
{s}_{\max
}^{T}(\bolds{\nu}-\bolds{\mu})$ that corresponds to the largest
``distance'' between the mean vectors,
that is,
the direction which best separates their distributions. The resulting
test,\vspace*{1pt}
known as Roy's largest root test, is given by the largest eigenvalue of
$\mathbf{BS}^{-1}$ where $\mathbf{B}$ is the matrix of between groups
(or populations) sums of squares and cross products, and $\mathbf
{S}$ is
the usual unbiased estimator of $\bolds{\Sigma}$. In the special case
when there are only two populations, this test statistic is identical to
Hotelling's $T^{2}$ statistic. From the simultaneous confidence intervals
point of view, the critical values derived from the null distribution
of this
statistic can be used for constructing Scheffe's simultaneous confidence
intervals for all possible linear combinations of the difference
$(\bolds{\mu}-\bolds{\nu})$. Further note that the estimated
direction corresponding to Roy's largest root test is $\mathbf{S}%
^{-1}(\overline{\mathbf{Y}}-\overline{\mathbf{X}})$ where
$\overline{\mathbf{X}}$ and $\overline{\mathbf{Y}}$ are the respective
sample means. 

Our objective is to extend and generalize the classical multivariate
method, described above, to nonnormal multivariate ordered data. Our
approach will be nonparametric. Recall\vspace*{1pt} that comparing MVNs
is done by considering the family of statistics
$T_{n,m}(\mathbf{s})=\mathbf{s}^{T}(\overline
{\mathbf{Y}}-\overline{\mathbf{X}})$ for all $\mathbf{s}\in $ $
\mathbb{R}^{p}$. In the case of nonnormal populations, the population
mean alone may not be enough to characterize the distribution. In such
cases, it may not be sufficient to compare the means of the populations
but one may have to compare entire distributions. One possible way of
doing so is by considering rank statistics. Suppose
$\mathbf{X}_{1},\ldots,\mathbf{X}_{n}$ and
$\mathbf{Y}_{1},\ldots,\mathbf{Y}_{m}$ are independent random samples
from two multivariate populations. Let
\[
R_{k} ( \mathbf{s} ) =\sum_{i=1}^{n}
\mathbb{I}_{ (
\mathbf{s}^{T}\mathbf{X}_{i}\leq\mathbf{s}^{T}\mathbf{X}
_{k} ) }+\sum_{j=1}^{m}
\mathbb{I}_{ ( \mathbf{s}^{T}%
\mathbf{Y}_{j}\leq\mathbf{s}^{T}\mathbf{X}_{k} ) }%
\]
be the rank of $\mathbf{s}^{T}\mathbf{X}_{k}$ in the combined sample
$\mathbf{s}^{T}\mathbf{X}_{1},\ldots,\mathbf
{s}^{T}\mathbf{X}%
_{n},\mathbf{s}_{1}^{T}\mathbf{Y}_{1},\ldots,\mathbf{s}%
^{T}\mathbf{Y}_{m}$. For fixed $\mathbf{s}\in
\mathbb{R}^{p}$ the distributions of $\mathbf{s}^{T}\mathbf{X}$ and
$\mathbf{s}^{T}\mathbf{Y}$ can be compared using a rank test. For
example, if we use $W_{n,m}(\mathbf{s})=\sum_{i=1}^{n}R_{i} (
\mathbf{s} ) $ our comparison is done in terms of Wilcoxon's rank
sum statistics. It is well known that rank tests are well suited for testing
for univariate stochastic order [cf. \citet{HajSidSen99}, \citet{Dav12}] where
the restrictions that $\mathbf{s}\in\mathcal{S}_{+}^{p-1}$ must be made.
Although any rank test can be used, the Mann--Whitney form of Wilcooxon's
(WMW) statistic is particularly attractive in this application.
Therefore in
the rest of this paper we develop estimation and testing procedures for the
linear stochastic order based on the family of statistics%
%
\begin{equation}\label{Psi-nm}
\Psi_{n,m}(\mathbf{s})=\frac{1}{nm}\sum
_{i=1}^{n}\sum_{j=1}^{m}%
\mathbb{I}_{ ( \mathbf{s}^{T}\mathbf{X}_{i}\leq\mathbf
{s}%
^{T}\mathbf{Y}_{j} ) },
\end{equation}
where $\mathbf{s}$ varies over $\mathcal{S}_{+}^{p-1}$. Note that
(\ref{Psi-nm}) unbiasedly estimates%
%
\begin{equation}\label{Psi}
\Psi(\mathbf{s})=\mathbb{P} \bigl( \mathbf{s}^{T}\mathbf{X}%
\leq\mathbf{s}^{T}\mathbf{Y} \bigr).%
\end{equation}
The following result is somewhat surprising.
%
\begin{proposition}
\label{Prop-Smax}Let $\mathbf{X}$ and $\mathbf{Y}$ be
independent MVNs
with means $\bolds{\mu}\leq$ $\bolds{\nu}$ and common variance
matrix $\bolds{\Sigma}$. Then Roy's maximal separating direction
$\bolds{\Sigma}^{-1}(\bolds{\nu}-\bolds{\mu})$ also maximizes
$\mathbb{P} ( \mathbf{s}^{T}\mathbf{X}\leq\mathbf{s}%
^{T}\mathbf{Y} )$.
\end{proposition}

Proposition~\ref{Prop-Smax} shows that the direction which separates the
means, in the sense of Roy, also maximizes (\ref{Psi}). Thus it provides
further support for choosing (\ref{Psi-nm}) as our test statistic. Note that
in general $\bolds{\Sigma}^{-1}(\bolds{\nu}-\bolds{\mu
})$ may
not belong to $\mathcal{S}_{+}^{p-1}$. Since we focus on the linear
statistical order, we restrict ourselves to $\mathbf{s}\in\mathcal
{S}%
_{+}^{p-1}$. Consequently we define $\mathbf{s}_{\max}:=\arg
\max_{\mathbf{s}\in\mathcal{S}_{+}^{p-1}}\Psi(\mathbf{s})$ and refer
to $\mathbf{s}_{\max}$ as the best separating direction. Further
note that
if $\mathbf{X}$ and $\mathbf{Y}$ are independent and continuous
and if
$\mathbf{X}\preceq_{\lst }\mathbf{Y}$, then $\Psi(\mathbf{s}%
)\geq1/2$ for all $\mathbf{s}\in\mathcal{S}_{+}^{p-1}$. This simply means
that $\mathbf{s}^{T}\mathbf{X}$ tends to be smaller than
$\mathbf{s}^{T}\mathbf{Y}$ more than $50\%$ of the time. Note that
probabilities of type (\ref{Psi}) were introduced by \citet{Pit37} and
further studied by \citet{Das85} for comparing estimators. Random variables
satisfying such a condition are said to be ordered by the precedence order
[\citet{ArcKvaSam02}].

Once $\mathbf{s}_{\max}$ is estimated we can plug it into (\ref
{Psi-nm}) to
get a test statistic. Hence our test may be viewed as a natural generalization
of Roy's largest root test from MVNs to arbitrary ordered distributions.
However, unlike Roy's method, which does not explicitly estimate
$\mathbf{s}_{\max}$, we do. On the other hand the proposed test
does not
require the computation of the inverse of the sample covariance matrix whereas
Roy's test and Hotteling's $T^{2}$ test require such computations.
Consequently, such tests cannot be used when $n<p$ whereas our test can be
used in all such instances.
%
\begin{remark}
In the above description $\mathbf{X}_{i}$ and $\mathbf{Y}_{j}$ are
independent for all $i$ and $j$ and therefore the probability $\mathbb
{P}%
( \mathbf{s}^{T}\mathbf{X}_{i}\leq\mathbf{s}^{T}%
\mathbf{Y}_{j} ) $ is independent of both $i$ and~$j$. However, in
many applications such as repeated measurement and crossover designs,
the data are a random sample of dependent pairs $ (
\mathbf{X}_{1},\mathbf{Y}_{1} ),\ldots, ( \mathbf{X}_{N},\mathbf
{Y}%
_{N} ) $ for which $\mathbf{Z}_{i}=\mathbf
{Y}_{i}-\mathbf{X}%
_{i}$ are i.i.d. For example, such a situation may arise when
$\mathbf{Y}%
_{i}=\mathbf{Y}_{i}^{\prime}+\bolds{\epsilon}_{i}$ and
$\mathbf{X}_{i}=\mathbf{X}_{i}^{\prime}+\bolds{\epsilon}_{i}$,
where $\bolds{\epsilon}_{i}$ are pair-specific random effects and the
RVs $\mathbf{Y}_{i}^{\prime}$ (as well as $\mathbf
{X}_{i}^{\prime}$)
are i.i.d. In this situation $\mathbb{P}(\mathbf{s}^{T}\mathbf
{X}%
_{i}\leq\mathbf{s}^{T}\mathbf{Y}_{i})$ is independent of $i$ and
$\mathbf{s}_{\max}$ is well defined. Moreover the objective function
analogous to (\ref{Psi-nm}) is%
%
\begin{equation}\label{Psi-N}
\Psi_{N}(\mathbf{s})=\frac{1}{N}\sum
_{i=1}^{N}\mathbb{I}_{ (
\mathbf{s}^{T}\mathbf{X}_{i}\leq\mathbf{s}^{T}\mathbf{Y}
_{i} ) }.
\end{equation}
In the following we consider both sampling designs which we refer to
as: (a)
independent samples and (b) paired or dependent samples. Results are
developed primarily for independent samples, but modification for paired
samples are mentioned as appropriate.
\end{remark}

\subsection{Estimating the best separating direction}\label{sec3.2}

Consider first the case of independent samples, that is, $\mathbf{X}
_{1},\ldots,\mathbf{X}_{n}$ and $\mathbf{Y}_{1},\ldots,\mathbf{Y}%
_{m}$ are random samples from the two populations. Rewrite (\ref
{Psi-nm}) as%
%
\begin{equation}\label{Psi-nm-Z}
\Psi_{n,m}(\mathbf{s})=\frac{1}{nm}\sum
_{i=1}^{n}\sum_{j=1}^{m}%
\mathbb{I}_{ ( \mathbf{s}^{T}\mathbf{Z}_{ij}\geq0 )},
\end{equation}
where $\mathbf{Z}_{ij}=\mathbf{Y}_{j}-\mathbf{X}_{i}$. The
maximizer of (\ref{Psi-nm-Z}) is denoted by $\widehat{\mathbf
{s}}_{\max}
$, that is,%
%
\begin{equation}\label{s-max-hat}
\widehat{\mathbf{s}}_{\max}=\arg\max_{\mathbf{s}\in\mathcal{S}%
_{+}^{p-1}}\Psi_{n,m}(
\mathbf{s}).%
\end{equation}
Finding (\ref{s-max-hat}) with $\mathbf{s}\in\mathcal{S}_{+}^{p-1}$
is a
nonsmooth optimization problem. Consider first the situation where
$p=2$. In
this case we maximize (\ref{Psi-nm-Z}) subject to $\mathbf{s}%
\in\mathcal{S}_{+}^{1}=\{
(s_{1},s_{2})\dvtx s_{1}^{2}+s_{2}^{2}=1,(s_{1},s_{2}%
)\geq0\}$. Geometrically $\mathcal{S}_{+}^{1}$ is a quarter circle spanning
the first quadrant. Now let $\mathbf{Z}=(Z_{1},Z_{2})$, and without any
loss of generality assume that $\llVert\mathbf{Z}\rrVert
=1$. We
examine the behavior of the function $\mathbb{I}_{ ( \mathbf{s}%
^{T}\mathbf{Z}\geq0 ) }$ as a function of $\mathbf{s}$. Clearly
if $\mathbf{Z}\geq0$, that is, if $Z_{1}\geq0,Z_{2}\geq0$, then for all
$\mathbf{s}\in\mathcal{S}_{+}^{1}$ we have $\mathbb{I}_{ (
\mathbf{s}^{T}\mathbf{Z}\geq0 ) }=1$. In other words any value
of $\mathbf{s}$ on the arc $\mathcal{S}_{+}^{1}$ maximizes $\mathbb
{I}%
_{ ( \mathbf{s}^{T}\mathbf{Z}\geq0 ) }$. Similarly if
$\mathbf{Z}<0$ then for all $\mathbf{s}\in\mathcal{S}_{+}^{1}$
we have
$\mathbb{I}_{ ( \mathbf{s}^{T}\mathbf{Z}\geq0 ) }=0$ and
again the entire arc $\mathcal{S}_{+}^{1}$ maximizes $\mathbb{I}_{ (
\mathbf{s}^{T}\mathbf{Z}\geq0 ) }$. Now let $Z_{1}\geq0$ and
$Z_{2}<0$. It follows that $\mathbb{I}_{ ( \mathbf{s}^{T}%
\mathbf{Z}\geq0 ) }=1$ provided $\cos(\mathbf{s}^{T}%
\mathbf{Z})\geq0$. Thus $\mathbb{I}_{ ( \mathbf{s}^{T}%
\mathbf{Z}\geq0 ) }=1$ for all $\mathbf{s}$ on the arc
$[0,\theta]$ for some $\theta$. If $Z_{1}<0$ and let $Z_{2}\geq0$ the
situation is reversed and $\mathbb{I}_{ ( \mathbf{s}^{T}%
\mathbf{Z}\geq0 ) }=1$ for all angles $\mathbf{s}$ on the arc
$[\theta,\pi/2]$. The value of $\theta$ is given by (\ref{Z->theta}).
In other
words each $\mathbf{Z}$ is mapped to an arc on $\mathcal
{S}_{+}^{1}$ as
described above. Now, the function (\ref{Psi-nm-Z}) simply counts the number
of arcs covering each $\mathbf{s}\in\mathcal{S}_{+}^{1}$. The
maximizer of
(\ref{Psi-nm-Z}) lies in the region where the maximum number of arcs overlap.
Clearly this implies that the maximizer of (\ref{Psi-nm-Z}) is not
unique. A
quick way to find the maximizer is the following:
%
\begin{algorithm}
\label{Alg-Max2}Let $M$ denote the number of $\mathbf{Z}_{ij}$'s which
belong to the second or fourth quadrant. Map%
%
\begin{equation}\label{Z->theta}
\mathbf{Z}_{ij}\mapsto\theta_{ij}= \cases{ \pi/2-
\cos^{-1}(Z_{ij,1}), &\quad if $Z_{ij,1}\geq0,Z_{ij,2}<0$,
\cr
\cos^{-1}(Z_{ij,1})-\pi/2, &\quad if $Z_{ij,1}<0,Z_{ij,2}
\geq0$.} %
\end{equation}
Relabel and order the resulting angles as $\theta_{[1]}<\cdots
<\theta_{[ M]}$. Also define $\theta_{[0]}=0$ and $\theta
_{[
M+1]}=\pi/2$. Evaluate $\Psi_{n,m}(\mathbf{s}_{[i]})$ $i=1,\ldots,M$ where
$\mathbf{s}_{[i],1}=\cos(\theta_{[ i]})$ and $\mathbf{s}%
_{[i],2}=\sin(\theta_{[ i]})$. If a maximum is attained at
$\theta_{[ j]}$, then any value in $[\theta_{[ j-1]},\theta
_{[ j]}]$ or $[\theta_{[ j]},\theta_{[ j+1]}]$ maximizes
(\ref{Psi-nm}).
\end{algorithm}

In light of the above discussion we can be easily prove the following:
%
\begin{proposition}
\label{Prop-Max2}For $p=2$ Algorithm~\ref{Alg-Max2} maximizes (\ref{Psi-nm}).
\end{proposition}

In the general case, that is, for $p\geq3$, each observation
$\mathbf{Z}%
_{ij}$ is associated with a ``slice'' of $\mathcal{S}_{+}^{p-1}$. The
boundaries of the slice are the intersection of $\mathcal{S}_{+}^{p-1}$
and some half-plane. Note that when $p=2$ the slices are arcs. The
shape of the slice depends on the quadrant to which $\mathbf{Z}_{ij}$
belongs. The maximizer of (\ref{Psi-nm}) is again the value of
$\mathbf{s}$ which belongs to the largest number of slices. Although the
geometry of the resulting optimization problem is easy to understand,
we have not been able to devise a simple algorithm, which scales with
$p$, based on the ideas above. However, we have found that
(\ref{s-max-hat}) can be obtained by converting the data into polar
coordinates and then using the Nelder--Mead algorithm which does not
require the objective function to be differentiable. We emphasize that
this maximization process results in a single maximizer of
(\ref{Psi-nm}) and we do not attempt to find the entire set of
maximizers. For completness we note that there are methods for
optimizing (\ref{Psi-nm}) specifically designed for nonsmooth
problems. For more details see \citet{PriReaRob06} and
\citet{AudBecLeD08} and the references therein for both algorithms
and convergence results.
%
\begin{remark}
It is clear that the estimation procedure for paired samples is the
same as
for independent samples.
\end{remark}

\subsection{Large sample behavior}\label{sec3.3}

We find three different asymptotic regimes for $\widehat{\mathbf{s}}
_{\max}$ depending on distributional assumptions and the sampling scheme
(paired versus independent samples).

Note that the parameter space is the unit sphere not the usual Euclidian
space. There are several ways of dealing with this irregularity, one of which
is to re-express the last coordinate of $\mathbf{s}\in\mathcal{S}%
_{+}^{p-1}$ as $\mathbf{s}_{p}=\sqrt{1-s_{1}^{2}-\cdots
-s_{p-1}^{2}}$ and
consider the parameter space $\{\mathbf{s}\geq\mathbf{0}\dvtx s_{1}%
^{2}+\cdots+s_{p-1}^{2}\leq1\}$ which is a compact subset of~$
\mathbb{R}^{p-1}$. Clearly these parameterizations are equivalent, and
without any
further ambiguity we will denote them both by $\mathcal{S}_{+}^{p-1}$.
Thus in
the proofs below both views of $\mathcal{S}_{+}^{p-1}$ are used
interchangeably as convenient.

We begin our discussion with independent samples assuming continuous
distributions for both $\mathbf{X}$ and $\mathbf{Y.}$
%
\begin{theorem}
\label{Them-LST1}Let $\mathbf{X}$ and $\mathbf{Y}$ have continuously
differentiable densities. If $\Psi(\mathbf{s})$ is uniquely
maximized by
$\mathbf{s}_{\max}\in\mathrm{interior}(\mathcal{S}_{+}^{p-1})$. Then
$\widehat{\mathbf{s}}_{\max}$, the maximizer of (\ref{Psi-nm}), is
strongly consistent, that is, $\widehat{\mathbf{s}}_{\max}\stackrel
{\mathit{a.s.}}{\rightarrow}\mathbf{s}_{\max}$. Furthermore $\widehat
{\mathbf{s}}_{\max}=\mathbf{s}_{\max}+O_{p}(N^{-1/2})$ where
$N=n+m. $
Finally, if $n/ ( n+m ) \rightarrow\lambda\in( 0,1
), $
then%
\[
N^{1/2}(\widehat{\mathbf{s}}_{\max}-\mathbf{s}_{\max
})
\Rightarrow N ( 0,\bolds{\Sigma} ),
\]
where the matrix $\bolds{\Sigma}$ is defined in the body of the proof.
\end{theorem}

Although (\ref{Psi-nm}) is not continuous (nor differentiable) its
$U$-statistic structure guarantees that it is ``almost'' so [i.e., it is
continuous up to an $o_{p}(1/N)$ term], and therefore its maximizer converges
at a $\sqrt{N}$ rate to a normal limit [\citet{She93}]. We also note
that it
is difficult to estimate the asymptotic variance $\bolds{\Sigma}$
directly since it depends on unknown functions ($\nabla\psi_{j}$ and
$\nabla^{2}\psi_{j}$ for $j=1,2$
are defined in the body of the proof). Nevertheless bootstrap variance estimates are
easily derived.
%
\begin{remark}
Note that if either $\mathbf{X}$ or $\mathbf{Y}$ are
continuous RVs,
then $\Psi(\mathbf{s})$ is continuous. This is a necessary
condition for
the uniqueness of $\mathbf{s}_{\max}$.\vadjust{\goodbreak}
\end{remark}

We have not been able to find general condition(s) for a unique
maximizer for
$\Psi(\mathbf{s})$, although important sufficient conditions can be found.
For example:

\begin{proposition}
\label{Prop-Unique}If $\mathbf{Z}=\mathbf{Y}-\mathbf{X}$
and there
exist $\bolds{\delta}=\bolds{\nu}-\bolds{\mu}\geq
\mathbf{0}$ and $\bolds{\Sigma}$ so the distribution of
\[
\frac{\mathbf{s}^{T}\mathbf{Z}-\mathbf{s}^{T}\bolds{\delta}}%
{\sqrt{\mathbf{s}^{T}\bolds{\Sigma}
\mathbf{s}}}%
\]
is independent of $\mathbf{s}$, then the maximizer of $\Psi
(\mathbf{s}%
)$ is unique.
\end{proposition}

The condition above is satisfied by location scale families, and it may be
convenient to think of $\bolds{\delta}$ and $\bolds{\Sigma}$
as the
location and scale parameters for $\mathbf{Z}$. In general, however,
$\Psi(\mathbf{s})$ may not have a unique maximum nor be continuous. For
example, if both $\mathbf{X}$ and $\mathbf{Y}$ are discrete
RVs, then
$\Psi(\mathbf{s})$ is a step function. In such situations
$\mathbf{s}%
_{\max}$ is set valued, and we may denote it by $\mathbf{S}_{\max
}$. As we
have seen earlier the maximizer of $\Psi_{n,m}(\mathbf{s})$ is
always set
valued (typically, however, we find only one maximizer). Consider, for example,
the case where $\mathbb{P}(\mathbf{X}=(-1,-1))=\mathbb
{P}(\mathbf{X}%
=(1,1))=1/2$ and let $\mathbb{P}(\mathbf{Y}=(-1,-1))=1/2-\varepsilon
$, and
$\mathbb{P}(\mathbf{Y}=(1,1))=1/2+\varepsilon$ for some $\varepsilon>0$.
It is clear that $\mathbf{X}\prec_{\st }\mathbf{Y}$. Further note that
$\mathbf{Z}\in\{ (2,2),(0,0),(-2,-2) \}, $ and it follows that
$\Psi(\mathbf{s})$ is constant on $\mathcal{S}_{+}^{1}$ which
implies that
$\mathbf{S}_{\max}$ coincides with $\mathcal{S}_{+}^{1}$. Similarly
$\Psi_{n,m}(\mathbf{s})$ is constant on $\mathcal{S}_{+}^{1}$ and
therefore $\widehat{\mathbf{s}}_{\max}\in\mathbf{S}_{\max}$ for all
$n,m$. This means that consistency is guaranteed and the limiting distribution
is degenerate. More generally, we have:
%
\begin{theorem}
\label{Them-LST2}If $\mathbf{X}$ and $\mathbf{Y}$ have discrete
distributions with finite support and $\widehat{\mathbf{s}}_{\max}$
is a
maximizer of (\ref{Psi-nm}), then%
\[
\mathbb{P}(\widehat{\mathbf{s}}_{\max}\notin\mathbf{S}_{\max
})\leq
C_{1}\exp( -C_{2}N )
\]
for some positive constants $C_{1}$ and $C_{2}$.
\end{theorem}

Theorem~\ref{Them-LST2} shows that the probability that a maximizer of
$\Psi_{n,m}$ is not in $\mathbf{S}_{\max}$ is exponentially small
when the
underlying distributions of $\mathbf{X}$ and $\mathbf{Y}$ are
discrete. Hence $\widehat{\mathbf{s}}_{\max}$ is consistent and converges
exponentially fast. In fact the proof of Theorem~\ref{Them-LST2}
implies that
$\widehat{\mathbf{S}}_{\max}\rightarrow\mathbf{S}_{\max}$; that
is, the
set of maximizers of $\Psi_{n,m} ( \cdot) $ converges to the set
of maximizers of $\Psi( \cdot)$; that is, $\rho_{\mathrm{H}}(\widehat
{\mathbf{S}}_{\max},\mathbf{S}_{\max})\rightarrow0$ where $\rho_{\mathrm{H}}$
is the Hausdroff metric defined on compact sets. A careful reading of the
proof shows that Theorem~\ref{Them-LST2} also holds under paired samples.

Finally, we consider the case of continuous RVs under paired samples. Then
under the conditions of Theorem~\ref{Them-LST1} and provided the
density of
$\mathbf{Z}=\mathbf{Y}-\mathbf{X}$ is bounded we have:
%
\begin{theorem}
\label{Them-LST3}Under the above mentioned conditions $\widehat
{\mathbf{s}}_{\max}$, the maximizer of (\ref{Psi-N}), is strongly
consistent, that is, $\widehat{\mathbf{s}}_{\max}\stackrel
{\mathit{a.s.}}{\rightarrow
}\mathbf{s}_{\max}$,\vadjust{\goodbreak} converges at a cube root rate, that is,
$\widehat
{\mathbf{s}}_{\max}=\mathbf{s}_{\max}+O_{p} ( N^{-1/3}
) $,
and%
\[
N^{1/3}(\widehat{\mathbf{s}}_{\max}-\mathbf{s}_{\max
})
\Rightarrow\mathbf{W},
\]
where $\mathbf{W}$ has the distribution of the almost surely unique
maximizer of the process $\mathbf{s}\longmapsto-[Q ( \mathbf
{s}%
) +\mathbb{W}(\mathbf{s})]$ on $\mathcal{S}_{+}^{p-1}$ where
$Q ( \mathbf{s} ) $ is a quadratic function and $\mathbb{W}
(\mathbf{s})$ is a zero mean Gaussian process described in the body
of the proof.
\end{theorem}

Theorem~\ref{Them-LST3} shows that in paired samples
$\widehat{\mathbf{s}}_{\max}$ is consistent, but in contrast with Theorem
\ref{Them-LST1} it converges at a cube-root rate to a nonnormal
limit. The
cube root rate is due to the discontinuous nature of the objective function
(\ref{Psi-N}). General results dealing with this kind of asymptotics for
independent observations are given by \citet{KimPol90}. The main
difference between Theorems~\ref{Them-LST1} and~\ref{Them-LST3} is
that the
objective function (\ref{Psi-nm}) is smoothed by its $U$-statistic structure
while (\ref{Psi-N}) is not.

\subsection{\texorpdfstring{A confidence set for $\widehat{\mathbf{s}}_{\max}$}{A confidence set for s max}}\label{sec3.4}

Since the parameter space is the surface of a unit sphere it is natural to
define the $(1-\alpha)\times100\%$ confidence set for $\mathbf
{s}_{\max}$
centered at $\widehat{\mathbf{s}}_{\max}$ by
\[
\bigl\{\mathbf{s}\in\mathcal{S}_{+}^{p-1}\dvtx \widehat{
\mathbf{s}}_{\max} 
^{T}\mathbf{s}\leq C_{\alpha,N}
\bigr\},
\]
where $C_{\alpha,N}$ satisfies $\mathbb{P}(\widehat{\mathbf
{s}}_{\max}%
^{T}\mathbf{s}\leq C_{\alpha,N})=1-\alpha$. For more details see
\citet{FisHal89} or \citet{PedCha96}. Hence the confidence
set is the set of all $\mathbf{s}\in\mathcal{S}_{+}^{p-1}$ which have a
small angle with $\widehat{\mathbf{s}}_{\max}$. In theory one may appeal
to Theorem~\ref{Them-LST1} to derive the critical value for any
$\alpha\in(0,1)$. However the limit law in Theorem~\ref{Them-LST1}
requires knowledge of unknown parameters and functions. For this
reason, we explore the bootstrap for estimating $C_{\alpha,N}$.
%
\begin{remark}
Since in the case of paired samples, the estimator converges at cube
root rate rather than the square root rate, the standard bootstrap
methodology may yield inaccurate coverage probabilities; see
\citet{AbrHua05} and \citet{SenBanWoo10}. For this reason we
recommend the ``M~out of N'' bootstrap methodology. For further
discussion on the ``M out of N'' bootstrap methodology one may refer to
\citet{Lee99}, \citet{DelRodWol01}, \citet{BicSak08}.
\end{remark}

\subsection{Testing for order}\label{sec3.5}

Consider first the case of independent samples where interest is in testing
the hypothesis
%
\begin{equation}\label{H0vsH1}
H_{0}\dvtx \mathbf{X}=_{\st }\mathbf{Y}\quad\mbox{versus}\quad
H_{1}\dvtx \mathbf{X}%
\prec_{\st }\mathbf{Y}.
\end{equation}
Thus (\ref{H0vsH1}) tests whether the distributions of $\mathbf{X}$ and
$\mathbf{Y}$ are equal or ordered (later on we briefly discuss testing
$H_{0}\dvtx \mathbf{X}\preceq_{\st }\mathbf{Y}$ versus
$H_{1}\dvtx\mathbf{X}\nprec_{\st }\mathbf{Y}$). In this section we propose
a new test for detecting an ordering among two multivariate
distributions based on the maximal separating direction. The test is
based on the following observation:
%
\begin{theorem}
\label{Thm-TestStat}Let $\mathbf{X}$ and $\mathbf{Y}$ be independent
and continuous RVs. If \mbox{$\mathbf{X}=_{\st }\mathbf{Y}$}, then
$\mathbb{P} ( \mathbf{s}^{T}\mathbf{X}\leq\mathbf{s}%
^{T}\mathbf{Y} ) =1/2$ for all $\mathbf{s}\in\mathcal{S}%
_{+}^{p-1}$, and if both \textup{(i)} $\mathbf{X}\preceq_{\st }\mathbf{Y}$ and
\textup{(ii)} $\mathbb{P} (
\mathbf{s}^{T}\mathbf{X}\leq\mathbf{s}^{T}\mathbf
{Y} )
>1/2$ for some $\mathbf{s}\in\mathcal{S}_{+}^{p-1}$ hold, then
$\mathbf{X}\prec_{\st }\mathbf{Y}$.
\end{theorem}

Theorem~\ref{Thm-TestStat} says that if it is known a priori that
$\{\mathbf{X}\preceq_{\st }\mathbf{Y}\}=\{\mathbf{X}=_{\st }%
\mathbf{Y}\}\cup\{\mathbf{X}\prec_{\st }\mathbf{Y}\}$, that
is, the RVs
are either equal or ordered [which is exactly what (\ref{H0vsH1}) implies],
then a strict linear stochastic ordering implies a strict ordering by the
usual multivariate stochastic order. In particular under the
alternative there
must exist a direction $\mathbf{s}\in\mathcal{S}_{+}^{p-1}$ for which
$\mathbf{s}^{T}\mathbf{X}\prec_{\lst }\mathbf
{s}^{T}\mathbf{Y}%
$.
%
\begin{remark}
\label{RmTestPhilosophy-I}The assumption that $\mathbf{X}\preceq
_{\st }\mathbf{Y}$ is natural in applications such as environmental sciences
where high exposures are associated with increased risk. Nevertheless
if the
assumption that $\mathbf{X}\preceq_{\st }\mathbf{Y}$ is not warranted
then the alternative hypothesis formulated in terms of the linear stochastic
order actually tests whether there exists a $\mathbf{s}\in\mathcal
{S}%
_{+}^{p-1}$ for which $\mathbb{P} ( \mathbf{s}^{T}\mathbf{X}
\leq\mathbf{s}^{T}\mathbf{Y} ) >1/2$. This amounts to a
precedence (or Pitman) ordering between $\mathbf{s}^{T}\mathbf
{X}$ and
$\mathbf{s}^{T}\mathbf{Y}$.
\end{remark}
%
\begin{remark}
\label{RmTestPhilosophy-II}In the proof of Theorem~\ref{Thm-TestStat}
we use
the fact that given that $\mathbf{X}\preceq_{\st }\mathbf{Y}$ we have
$\mathbf{X}\prec_{\st }\mathbf{Y}$ provided $X_{i}\prec_{\st }Y_{i}$ for
some $1\leq i\leq p$. Note that if $X_{i}\prec_{\st }Y_{i}$, then $\mathbb
{E}(X_{i})<\mathbb{E}(Y_{i})$. Thus it is possible to test (\ref{H0vsH1}) by
comparing means (or any other monotone function of the data). Although
such a
test will be consistent it may lack power because tests based on means are
often far from optimal when the data is not normally distributed. The WMW
procedure, however, is known to have high power for a broad collection of
underlying distributions.
\end{remark}

Hence (\ref{H0vsH1}) can be reformulated in terms of the linear
stochastic. In
particular it justifies using the statistic%
%
\begin{equation}\label{Snm}
S_{n,m}=N^{1/2}\bigl(\Psi_{n,m}(\widehat{
\mathbf{s}}_{\max})-1/2\bigr).%
\end{equation}
To the best of our knowledge this is the first general test for multivariate
ordered distributions. In practice we first estimate $\widehat
{\mathbf{s}%
}_{\max}$ and then define $\widehat{U}_{i}=\widehat{\mathbf
{s}}_{\max}%
^{T}\mathbf{X}_{i}$ and $\widehat{V}_{j}=\widehat{\mathbf
{s}}_{\max
}^{T}\mathbf{Y}_{j}$ where $i=1,\ldots,n$ and $j=1,\ldots,m$. Hence
(\ref{Snm}) is nothing but a WMW test based on the $\widehat{U}{}^{\prime
}$'s and
$\widehat{V}$'s. It is also a Kolmogorov--Smirnov type test.

The large sample distribution of (\ref{Snm}) is given in the following.
%
\begin{theorem}
\label{Thm-TestLimit}Suppose the null (\ref{H0vsH1}) holds. Let
$n,m\rightarrow\infty$ and $n/(n+m)\rightarrow\lambda\in(0,1)$. Then
\[
S_{n,m}\Rightarrow S=\sup_{\mathbf{s}\in\mathcal{S}_{+}^{p-1}}%
\mathbb{G}(
\mathbf{s}),
\]
where $\mathbb{G}(\mathbf{s})$ is a zero mean Gaussian process with
covariance function given by (\ref{Cu,v-H0}).\vadjust{\goodbreak}
\end{theorem}
%
\begin{remark}
Since $\widehat{\mathbf{s}}_{\max}\stackrel{\mathrm{a.s.}}{\rightarrow
}\mathbf{s}_{\max}$ by Slutzky's theorem the power of test (\ref{Snm})
converges to the power of a WMW test comparing the samples $(\mathbf
{s}%
_{\max}^{T}\mathbf{X}_{1},\ldots,\mathbf{s}_{\max
}^{T}\mathbf{X}%
_{n})$ and $(\mathbf{s}_{\max}^{T}\mathbf{Y}_{1},\ldots,\mathbf{s}%
_{\max}^{T}\mathbf{Y}_{m})$. The ``synthetic'' test, assuming that
$\mathbf{s}_{\max}$ is known, serves as a gold standard as verified
by our
simulation study.
\end{remark}
%
\begin{remark}
Furthermore, the power of the test under local alternatives, that is, when
$\mathbf{Y}=_{\st }\mathbf{X}+N^{-1/2}\bolds{\delta}$ and
$N\rightarrow\infty$ is bounded by the power of the WMW test comparing the
distributions of $\mathbf{s}_{\max}^{T}\mathbf{X}$ and
$\mathbf{s}%
_{\max}^{T}\mathbf{Y}=\mathbf{s}_{\max}^{T}\mathbf{X}+
N^{-1/2}\mathbf{s}_{\max}^{T}\bolds{\delta}$.
\end{remark}

Alternatives to the ``sup'' statistic (\ref{Snm}) are the
``integrated''
statistics%
%
\begin{eqnarray}\label{Int-Tests}
I_{n,m}&=&\int_{\mathbf{s}\in\mathcal{S}_{+}^{p-1}}\bigl[N^{1/2}\bigl(
\Psi_{n,m}(\mathbf{s})-1/2\bigr)\bigr]\,d\mathbf{s}
\quad\mbox{and}\nonumber\\[-8pt]\\[-8pt]
I_{n,m}^{+}&=&\int_{\mathbf{s}\in\mathcal{S}_{+}^{p-1}}
\bigl[N^{1/2}\bigl(\Psi_{n,m}%
(\mathbf{s})-1/2\bigr)
\bigr]_{+}\,d\mathbf{s},\nonumber
\end{eqnarray}
where $ [ x ]_{+}=\max( 0,x )$. It is clear that
$I_{n,m}\Rightarrow N ( 0,\sigma^{2} ) $ where
\[
\sigma^{2}=\int
_{\mathbf{u}\in\mathcal{S}_{+}^{p-1}}\int_{\mathbf{v}%
\in\mathcal{S}_{+}^{p-1}}C ( \mathbf{u},\mathbf{v} )
\,d\mathbf{u}\,d\mathbf{v}
\]
and $C ( \mathbf{u},\mathbf
{v}%
)$, the covariance function of $\mathbb{G}$, is given by
(\ref{Cu,v-H0}). Also%
\[
I_{n,m}^{+}\Rightarrow\int_{\mathbf{s}\in\mathcal{S}_{+}^{p-1}%
}\bigl[
\mathbb{G}(\mathbf{s})\bigr]_{+}\,d\mathbf{s}.
\]
This distribution does not have a closed form. The statistics $I_{n,m}$ and
$I_{n,m}^{+}$ have univariate analogues; cf. \citet{DavHer12}.
Finally, we have the following theorem:
%
\begin{theorem}
\label{Thm-Consistent&Monotone}The tests (\ref{Snm}) and (\ref
{Int-Tests}) are
consistent. Furthermore if $\mathbf{X}\preceq_{\lst }\mathbf{Y}%
\preceq_{\lst }\mathbf{Z}$, then all three tests for
$H_{0}\dvtx \mathbf{X}%
=_{\st }\mathbf{Z}$ versus $H_{1}\dvtx \mathbf{X}\prec_{\st }\mathbf{Z}$
are more powerful than the respective tests for $H_{0}\dvtx \mathbf{X}%
=_{\st }\mathbf{Y}$ versus $H_{1}\dvtx \mathbf{X}\prec_{\st }\mathbf{Y}$.
\end{theorem}

Theorem~\ref{Thm-Consistent&Monotone} shows that the tests are
consistent and
that their power function is ``monotone'' in the linear stochastic order.
%
\begin{remark}
Qualitatively similar results are obtainable in the paired sampling
case; the
only difference being the limiting process. For example, it easy to see that
the paired sample analogue of (\ref{Snm}) satisfies
\[
N^{1/2}\bigl(\Psi_{N}(\widehat{\mathbf{s}}_{\max})-1/2
\bigr)\Rightarrow\sup_{\mathbf{s}\in\mathcal{S}_{+}^{p-1}}\mathbb
{Q}(\mathbf{s}),
\]
where $\mathbb{Q}(\mathbf{s})$ is the empirical process on $\mathcal
{S}%
_{+}^{p-1}$ associated with (\ref{Psi-N}). Analogues of $I_{n,m}$ and
$I_{n,m}^{+}$ are similarly defined and analyzed. Tests for paired samples
may be similarly implemented using bootstrap or permutation methods.
\end{remark}

\section{Simulations}\label{sec4}

For simplicity of exposition, and motivated by the fact that the
example we
analyzed in this paper deals with independent samples, we limit our
simulations to the case of independent samples.

\subsection{\texorpdfstring{The distribution of $\widehat{\mathbf{s}}_{\max}$}{The distribution of s max}}\label{sec4.1}

\begin{sidewaysfigure}

\includegraphics{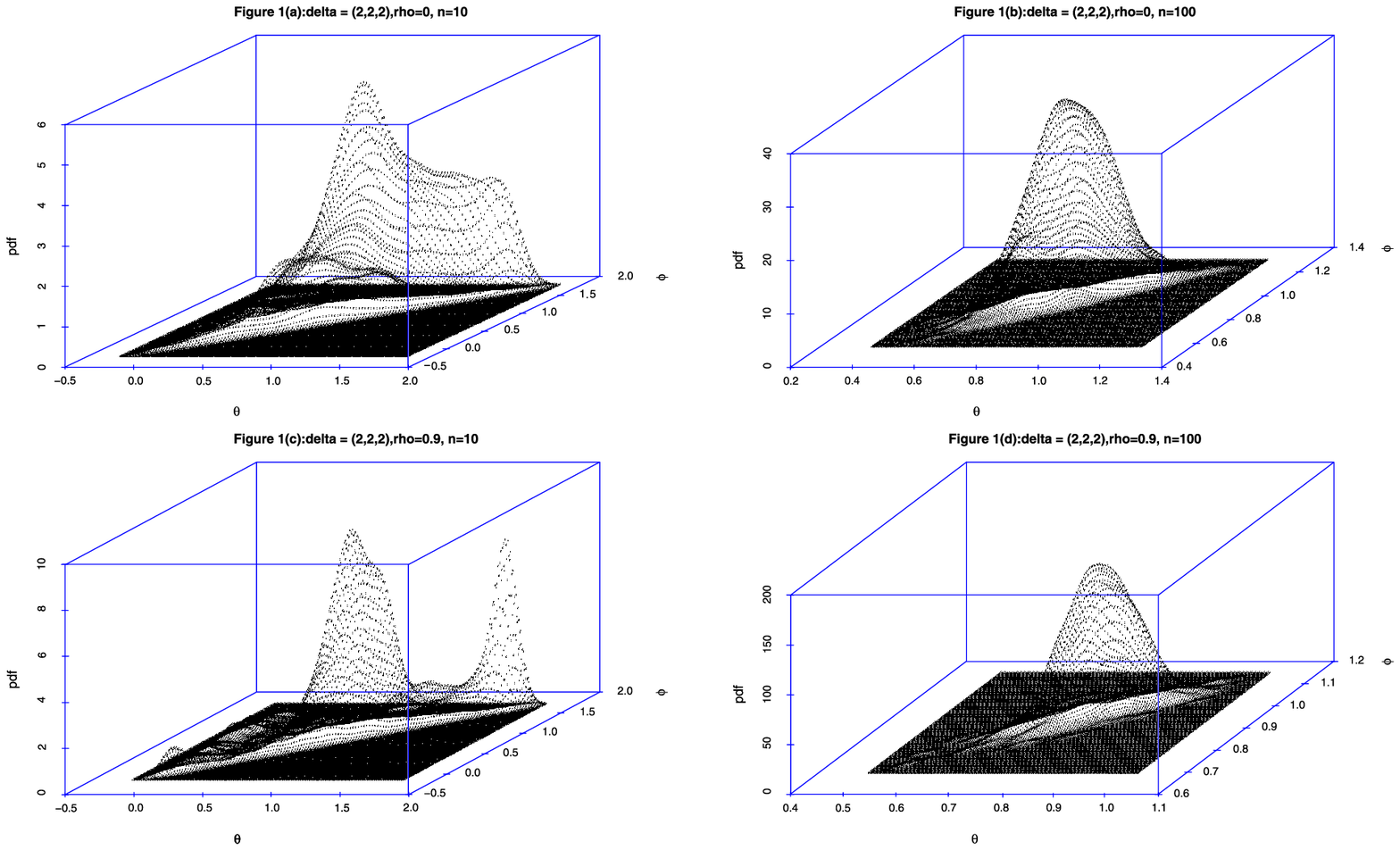}

\caption{Plot of (simulated) polar angles.}\label{figur1}
\end{sidewaysfigure}

\begin{sidewaysfigure}

\includegraphics{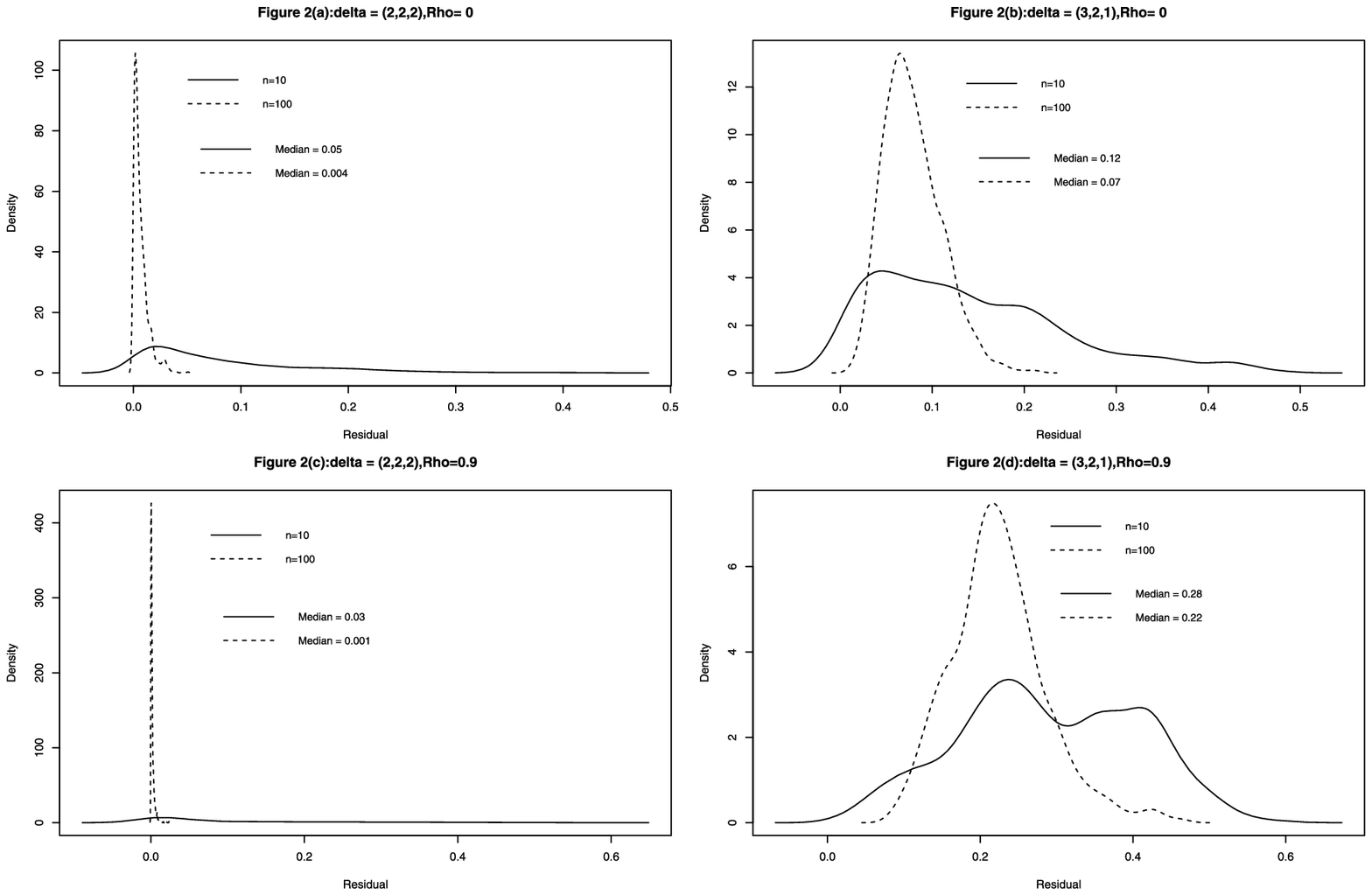}

\caption{Plot of (simulated) polar residuals.}\label{figur2}
\end{sidewaysfigure}

We start by investigating the distribution of $\widehat{\mathbf
{s}}_{\max
}$ by simulation. For simplicity we choose $p=3$ and generated
$\mathbf{X}%
_{i}$ ($i=1,\ldots,n$) distributed as $N_{3} (
\mathbf{0},\bolds{\Sigma} ) $ and $\mathbf{Y}_{j}$ ($j=1,\ldots,m$)
distributed as $N_{3} ( \bolds{\delta},\bolds{\Sigma} ) $
where $\bolds{\Sigma}= ( 1-\rho) \mathbf{I}%
+\rho\mathbf{J}$, $\mathbf{I}$ is the identity matrix and $\mathbf{J}$
is a square matrix of $1s$. We simulated $1000$ realizations of
$\widehat{\mathbf{s}}_{\max}$ for various sample sizes and correlation
coefficients. To get a visual description of the density of
$\widehat{\mathbf{s}}_{\max}$, we provide a pair of plots for each
configuration of $\rho$ and sample size $n$. In Figure~\ref{figur1} we
provide the joint density of the two-dimensional polar angles $ (
\theta,\phi) $ of $\widehat{\mathbf{s}}_{\max}$. There are four panels
in Figure~\ref{figur1}, corresponding to all combinations of
$\rho=0,0.9$ and $n=10,100$. The mean vector $\bolds{\delta}$ in this
plot was taken to be $\bolds {\delta }= ( 2,2,2 )^{T}$. In
Figure~\ref{figur2} we provide the density of the polar residual
defined by $1-\widehat{\mathbf{s}}^{T}_{\max}$
$\mathbf{s}%
_{\max}$. The four panels of Figure~\ref{figur2} correspond to all
combinations of $\rho=0,0.9$ and $n=10,100$ and two patterns of
$\bolds{\delta}$, namely, $(2,2,2)^{T}$ and $(3,2,1)^{T}$. We see from
Figure~\ref{figur1} that $\widehat {\mathbf{s}}_{\max}$ converges to a
unimodal, normal looking distribution as the sample size increases.
Interestingly, from Figure~\ref{figur2} we see that the concentration
of the distribution around the true parameter depends upon the values
of $\bolds{\delta}$ and $\rho$ (which together
determine~$\mathbf{s}_{\max}$). If the components of the underlying
random vector are exchangeable [e.g., $\bolds{\delta}= ( 2,2,2 )^{T}$],
the residuals tend to concentrate more closely around zero [Figure
\ref{figur2}(a) and (c)] compared to the case when they are not
exchangeable [Figure~\ref{figur2}(b) and
(d)].%

\subsection{Study design}\label{sec4.2}

The simulation study consists of three parts. In the first part we evaluate
the accuracy and precision of $\widehat{\mathbf{s}}_{\max}$ by estimating
its bias and mean squared error (MSE). In the second part we
investigate the
coverage probability of bootstrap confidence intervals. In the third
part we
estimate type I errors and powers of the proposed test $S_{n,m}$ as
well as
the integral tests $I_{n,m}$ and $I_{n,m}^{+}$.

To evaluate the bias and MSEs we generated
$\mathbf{X}_{1},\ldots,\mathbf{X}_{n}\sim N_{3} (
\mathbf{0},\bolds{\Sigma} ) $ and
$\mathbf{Y}_{1},\ldots,\mathbf{Y}_{m}\sim N_{3} (
\bolds{\delta},\bolds{\Sigma} ) $ where $n=m=20$ or $100$
observations. The common variance matrix is assumed to have intra-class
correlation structure, that is, $\bolds{\Sigma}= ( 1-\rho)
\mathbf{I}+\rho\mathbf{J}$ where $\mathbf{I}$ is the identity matrix and
$\mathbf{J}$ is a matrix of ones. Various patterns of the mean vectors
$\bolds{\delta}$ and correlation coefficient $\rho$ were considered as
described in Table~\ref{tabl1}.

We conducted extensive simulation studies to evaluate the performance
of the
bootstrap confidence intervals. In this paper we present a small sample
of our
study. We generated data from two 5-dimensional normal populations
with means
$\mathbf{0}$ and~$\bolds{\delta}$, respectively, and a common
covariance $\bolds{\Sigma}= ( 1-\rho) \mathbf{I}%
+\rho\mathbf{J}$. We considered 5 patterns of $\rho$ and 2 patterns of
sample sizes ($n=m=20$ and $40$). The nominal coverage probability was $0.95.
$ Results are summarized in Table~\ref{tabl2}.

%
\begin{table}
\tablewidth=270pt
\caption{Bias and MSE of $\widehat{\mathbf{s}}_{\max}$}\label{tabl1}
\begin{tabular*}{\tablewidth}{@{\extracolsep{\fill}}ld{2.2}d{1.5}c@{}}
\hline
$\bolds{\delta}$ & \multicolumn{1}{c}{$\bolds{\rho}$} & \multicolumn{1}{c}{\textbf{Bias}} & \textbf{MSE}\\
\hline
\multicolumn{4}{@{}c@{}}{$n=m=20$}\\[4pt]
$ ( 1,1,1 ) $ & -0.25 & 0.001 & 0.072\\
$ ( 1,1,1 ) $ & 0 & 0.004 & 0.129\\
$ ( 1,1,1 ) $ & 0.25 & 0.009 & 0.187\\
$ ( 1,1,1 ) $ & 0.50 & 0.012 & 0.216\\
$ ( 1,1,1 ) $ & 0.90 & 0.010 & 0.203\\[4pt]
$ ( 3,2,1 ) $ & -0.25 & 0.018 & 0.090\\
$ ( 3,2,1 ) $ & 0 & 0.001 & 0.066\\
$ ( 3,2,1 ) $ & 0.25 & 0.053 & 0.114\\
$ ( 3,2,1 ) $ & 0.50 & 0.060 & 0.113\\
$ ( 3,2,1 ) $ & 0.90 & 0.112 & 0.170\\[4pt]
\multicolumn{4}{@{}c@{}}{$n=m=100$}\\[4pt]
$ ( 1,1,1 ) $ & -0.25 & 0.00009 & 0.014\\
$ ( 1,1,1 ) $ & 0 & 0.00021 & 0.027\\
$ ( 1,1,1 ) $ & 0.25 & 0.00045 & 0.041\\
$ ( 1,1,1 ) $ & 0.50 & 0.00079 & 0.056\\
$ ( 1,1,1 ) $ & 0.90 & 0.00050 & 0.044\\[4pt]
$ ( 3,2,1 ) $ & -0.25 & 0.02400 & 0.039\\
$ ( 3,2,1 ) $ & 0 & 0.00004 & 0.012\\
$ ( 3,2,1 ) $ & 0.25 & 0.05200 & 0.065\\
$ ( 3,2,1 ) $ & 0.50 & 0.06400 & 0.077\\
$ ( 3,2,1 ) $ & 0.90 & 0.14100 & 0.158\\
\hline
\end{tabular*}
\end{table}

\begin{table}
\tablewidth=260pt
\caption{Coverage probabilities for the bootstrap confidence
intervals for $p=5$ normal data. Pattern $i=1,2$ corresponds to
$\bolds{\delta}_{1}= (0.1,0.25,0.5,0.75,0.9 ) $ and
$\bolds{\delta}_{2}= (0.5,0.5,0.5,0.5,0.5 )$}\label{tabl2}
\begin{tabular*}{\tablewidth}{@{\extracolsep{\fill}}ld{2.2}cc@{}}
\hline
\multicolumn{2}{@{}c}{\textbf{Set up}} & \multicolumn{2}{c@{}}{\textbf{Coverage
probability}}\\[-4pt]
\multicolumn{2}{@{}c}{\hrulefill} & \multicolumn{2}{c@{}}{\hrulefill}
\\
\textbf{Pattern} & \multicolumn{1}{c}{$\bolds{\rho}$}
& \multicolumn{1}{c}{$\bolds{n=m=20}$} & \multicolumn{1}{c@{}}{$\bolds{n=m=40}$}\\
\hline
1 & -0.25 & 0.981 & 0.971\\
1 & 0 & 0.913 & 0.918\\
1 & 0.25 & 0.916 & 0.933\\
1 & 0.50 & 0.971 & 0.969\\
1 & 0.90 & 0.993 & 0.989\\
[4pt]
2 & -0.25 & 0.982 & 0.967\\
2 & 0 & 0.984 & 0.972\\
2 & 0.25 & 0.986 & 0.978\\
2 & 0.50 & 0.968 & 0.968\\
2 & 0.90 & 0.950 & 0.954\\
\hline
\end{tabular*}
\end{table}

The goal of the third part of our simulation study is to evaluate the
type I error and the power of the test (\ref{Snm}). To evaluate the
type I error three different baseline distributions for the two
populations $\mathbf{X}$ and $\mathbf{Y}$ were employed as follows: $(1)$
both distributed as $N ( \mathbf{0},\bolds{\Sigma} )
\mathbf{;}(2)$ both distributed as $\pi N (
\mathbf{0},\bolds{\Sigma} ) + ( 1-\pi) N (
\bolds{\delta},\bolds{\Sigma} ) $ with $\pi=0.2$ or
$\pi=0.8\mathbf{;}$ and $ ( 3 ) $ both distributed as
$\exp( \mathbf{Z} ) = ( \exp(
Z_{1} ),\ldots,\exp( Z_{p} ) ) $ where $\mathbf{Z}$
follows a $N ( \bolds{\delta},\bolds{\Sigma} )$. We refer to
this distribution as the multivariate lognormal distribution.
Throughout the variance matrix is assumed to have the intra-class
structure described above. Various patterns of the mean vectors
$\bolds{\delta}$ and correlation coefficient $\rho$ the dimension $p$
were considered as described in Table~\ref{tabl3}. Sample sizes of $n=m=15$ or
$25$ are reported.

\begin{table}
\tablewidth=270pt
\caption{Type I errors for the proposed procedure with\break nominal
level $\alpha=0.05$. Three types of distributions are considered: MVNs,
MV-LogN
(multivariate lognormal) and Mix-MVN (mixtures of MVNs)}\label{tabl3}
\begin{tabular*}{\tablewidth}{@{\extracolsep{\fill}}lcd{2.2}cc@{}}
\hline
\multicolumn{3}{@{}c}{\textbf{Set up}}
& \multicolumn{2}{c@{}}{\textbf{Type I error}}\\
\hline
\textbf{Distribution}
& \multicolumn{1}{c}{$\bolds{p}$}
& \multicolumn{1}{c}{$\bolds{\rho}$}
& \multicolumn{1}{c}{$\bolds{n=m=15}$}
& \multicolumn{1}{c@{}}{$\bolds{n=m=25}$}\\
\hline
MVNs & 3 & -0.25 & 0.041 & 0.037\\[-0.5pt]
MVNs & 3 & 0.00 & 0.023 & 0.044\\[-0.5pt]
MVNs & 3 & 0.25 & 0.037 & 0.033\\[-0.5pt]
MVNs & 3 & 0.50 & 0.027 & 0.032\\[-0.5pt]
MVNs & 3 & 0.90 & 0.031 & 0.036\\
[4pt]
MVNs & 5 & -0.25 & 0.035 & 0.035\\[-0.5pt]
MVNs & 5 & 0.00 & 0.040 & 0.041\\[-0.5pt]
MVNs & 5 & 0.25 & 0.045 & 0.032\\[-0.5pt]
MVNs & 5 & 0.50 & 0.038 & 0.043\\[-0.5pt]
MVNs & 5 & 0.90 & 0.044 & 0.031\\
[4pt]
MV-LogN & 3 & -0.25 & 0.025 & 0.040\\[-0.5pt]
MV-LogN & 3 & 0.00 & 0.038 & 0.049\\[-0.5pt]
MV-LogN & 3 & 0.25 & 0.025 & 0.027\\[-0.5pt]
MV-LogN & 3 & 0.50 & 0.028 & 0.037\\[-0.5pt]
MV-LogN & 3 & 0.90 & 0.026 & 0.034\\
[4pt]
MV-LogN & 5 & -0.25 & 0.026 & 0.039\\[-0.5pt]
MV-LogN & 5 & 0.00 & 0.035 & 0.018\\[-0.5pt]
MV-LogN & 5 & 0.25 & 0.039 & 0.039\\[-0.5pt]
MV-LogN & 5 & 0.50 & 0.036 & 0.046\\[-0.5pt]
MV-LogN & 5 & 0.90 & 0.034 & 0.042\\
[4pt]
Mix-MVNs & 3 & -0.25 & 0.032 & 0.040\\[-0.5pt]
Mix-MVNs & 3 & 0.00 & 0.038 & 0.028\\[-0.5pt]
Mix-MVNs & 3 & 0.25 & 0.039 & 0.032\\[-0.5pt]
Mix-MVNs & 3 & 0.50 & 0.036 & 0.035\\[-0.5pt]
Mix-MVNs & 3 & 0.90 & 0.041 & 0.028\\
[4pt]
Mix-MVNs & 5 & -0.25 & 0.042 & 0.035\\[-0.5pt]
Mix-MVNs & 5 & 0.00 & 0.040 & 0.031\\[-0.5pt]
Mix-MVNs & 5 & -0.25 & 0.041 & 0.028\\[-0.5pt]
Mix-MVNs & 5 & 0.50 & 0.034 & 0.040\\[-0.5pt]
Mix-MVNs & 5 & 0.90 & 0.042 & 0.036\\
\hline
\end{tabular*} \vspace*{-3pt}
\end{table}

%
\begin{table}
\caption{Type I errors and power for some settings with
$p\geq n$.
Here $n=m=10$ and $\bolds{\delta}_{1}$ has components $1/2$ and
$\bolds{\delta}_{2}$ has components $i/p$}\label{tabl4}
\begin{tabular*}{\tablewidth}{@{\extracolsep{4in minus 4in}}lcd{1.3}@{\hspace*{-5pt}}
ccd{1.3}@{}}
\hline
\multicolumn{3}{@{}c}{\textbf{Type I error and power} $\bolds{p\!=\!10},\, \bolds{n\!=\!m\!=\!10}$} &
\multicolumn{3}{c@{}}{\hspace*{-15pt}\textbf{Type I error and power} $\bolds{p\!=\!20},\, \bolds{n\!=\!m\!=\!10}$}\\
\hline
$\bolds{\delta}$ & $\bolds{\rho}$ & \multicolumn{1}{c}{\textbf{Type I error}}
& $\bolds{\delta}$ & $\bolds{\rho}$
& \multicolumn{1}{c@{}}{\textbf{Type I error}}\\
\hline
$\mathbf{0}$ & 0.00 & 0.054 & $\mathbf{0}$ & 0.00 & 0.081\\
$\mathbf{0}$ & 0.25 & 0.051 & $\mathbf{0}$ & 0.25 & 0.050\\
$\mathbf{0}$ & 0.50 & 0.028 & $\mathbf{0}$ & 0.50 & 0.046\\
$\mathbf{0}$ & 0.90 & 0.038 & $\mathbf{0}$ & 0.90 &
0.048\\
\hline
\multicolumn{3}{@{}c}{\textbf{Power}} & \multicolumn{3}{c@{}}{\textbf{Power}}\\
\hline
$\bolds{\delta}_{1}$ & 0.00 & 0.83 & $\bolds{\delta}_{1}$ &
0.00 &
0.97\\
$\bolds{\delta}_{1}$ & 0.25 & 0.48 & $\bolds{\delta}_{1}$ &
0.25 &
0.53\\
$\bolds{\delta}_{1}$ & 0.50 & 0.26 & $\bolds{\delta}_{1}$ &
0.50 &
0.42\\
$\bolds{\delta}_{1}$ & 0.90 & 0.20 & $\bolds{\delta}_{1}$ &
0.90 &
0.22\\
[4pt]
$\bolds{\delta}_{2}$ & 0.00 & 0.98 & $\bolds{\delta}_{2}$ &
0.00 &
0.98\\
$\bolds{\delta}_{2}$ & 0.25 & 0.80 & $\bolds{\delta}_{2}$ &
0.25 &
0.59\\
$\bolds{\delta}_{2}$ & 0.50 & 0.67 & $\bolds{\delta}_{2}$ &
0.50 &
0.43\\
$\bolds{\delta}_{2}$ & 0.90 & 0.71 & $\bolds{\delta}_{2}$ &
0.90 &
0.40\\
\hline
\end{tabular*}
\end{table}

Power comparisons were carried out for data generated from $\mathbf
{X}%
_{1},\ldots,\mathbf{X}_{n}\sim N_{p} (
\mathbf{0},\bolds{\Sigma} ) $ and $\mathbf{Y}_{1},\ldots,\mathbf
{Y}%
_{m}\sim N_{p} ( \bolds{\delta},\bolds{\Sigma} ) $ where
$p=3$ or $5$ and a variety of patterns for $\bolds{\delta}$ as described
in Table~\ref{tabl4}. If Roy's maximal separating direction (cf. Proposition
\ref{Prop-Smax}) was known then a ``natural gold
standard'' would be the test based on $\Psi_{n,m}%
(\mathbf{s}_{\max})$. We shall refer to this test as the true maximal
direction (TMD) test. Clearly the TMD test cannot be used in practice
since it
involves the unknown direction $\mathbf{s}_{\max}$. Nevertheless
the TMD
test provides an upper bound for the power of the proposed test which
uses the
estimated direction. Hence we compute the efficiency of the proposed test
relative to TMD test. An additional test, referred to as the RMD test
is also
compared. The RMD test has the same form but uses Roy's maximal direction
given by $\mathbf{S}^{-1}(\overline{\mathbf{Y}}-\overline
{\mathbf{X}})$. As suggested by a reviewer we also evaluated the
power of
the two integral based tests, described in (\ref{Int-Tests}), which do not
require the determination of the best separating direction.

Additionally, in Table~\ref{tabl5} we evaluate the type I error and power of our
test when $\mathbf{X}_{1},\ldots,\mathbf{X}_{n}\sim N_{p} (
\mathbf{0},\bolds{\Sigma} ) $ and
$\mathbf{Y}_{1},\ldots,\mathbf{Y}_{m}\sim N_{p} (
\bolds{\delta},\bolds{\Sigma} ) $ and $n=m=p=10$ and $n=m=10$ and
$p=20$ (i.e., $p<n$ set up). Note that in neither of these cases the
standard Hotteling's $T^{2}$ (or Roy's largest root test) can be
computed whereas the proposed test can be calculated.

Simulation results reported in this paper are based on 1000 simulation runs.
Confidence sets are calculated using 1000 bootstrap samples. The bootstrap
critical values for estimating type I error were based on 500 bootstrap
samples. Since the results between 100 bootstrap samples and 500 bootstrap
samples did not differ by much, all powers were estimated using 100
bootstrap samples.\vspace*{-2pt}

\subsection{Simulation results}\label{sec4.3}

The Bias and MSEs for the patterns considered are summarized in Table~\ref{tabl1}.
It is clear that the bias decreases with the sample size as do
the\vadjust{\goodbreak}
MSEs. We
observe that the bias tends to be smaller under independence and negative
dependence compared with positive dependence. It also tends to be
smaller when
the data are exchangeable. Although results are not presented, we evaluated
squared bias and MSE for larger values of $p$ (e.g., $p=5,10$ and $20$)
and as
expected the total squared bias and total MSE increased with the dimension
$p$.

In Table~\ref{tabl2} we summarize the estimated coverage probabilities of the bootstrap
confidence intervals when $p=5$. Our simulation study suggests that the
proposed bootstrap methodology seems to perform better for larger
sample sizes
but rather poorly for smaller samples sizes.

Type I errors for different patterns considered in our simulation study are
summarized in Table~\ref{tabl3}.
Our simulation studies suggest that in every case the proposed
bootstrap based
test maintains the nominal level of 0.05. In general it is slightly
conservative. The performance of the test is not affected by the shape
of the
underlying distribution. This is not surprising, owing to the nonparametric
nature of the test. Furthermore, we evaluated the type I error of the proposed
bootstrap test for testing the null hypothesis (\ref{H0vsH1}) for $p$
as large
as $20$ with $n=m=10$ and discovered that the proposed test attains the
nominal level of 0.05 even $n\leq p$. See Table~\ref{tabl4}. As commented earlier
in the
paper, Hotelling's $T^{2}$ statistic cannot be applied here since the Wishart
matrix is singular in this case. However, the proposed method is still
applicable since the estimation of the best direction does not require the
inversion of a matrix.%

The power of tests (\ref{Snm}) and (\ref{Int-Tests}) for various patterns
considered in our simulation study are summarized in Table~\ref{tabl5}.

%
\begin{table}
\caption{Power comparisons of the two proposed test procedures
with type I error of 0.050. Here $\bolds{\delta}_{1}= (
0.1,0.5,0.9 )$, $\bolds{\delta}_{2}= (
0.1,0.25,0.5,0.75,0.9 )$, $\bolds{\delta}_{3}= (
0.5,0.5,0.5 ) $ and $\bolds{\delta}_{4}= (
0.5,0.5,0.5,0.5,0.5 ) $}\label{tabl5}
\begin{tabular*}{\tablewidth}{@{\extracolsep{\fill}}lcd{2.2}cccll@{}}
\hline
\multicolumn{8}{@{}c@{}}{\textbf{Power and RE \%} $\bolds{(n=m=15)}$}\\
\hline
\multicolumn{3}{@{}c}{\textbf{Set up}}
& \multicolumn{3}{c}{\textbf{Directional tests}} &
\multicolumn{2}{c@{}}{\textbf{Integral tests}}\\[-4pt]
\multicolumn{3}{@{}c}{\hrulefill}
& \multicolumn{3}{c}{\hrulefill} &
\multicolumn{2}{c@{}}{\hrulefill}\\
$\bolds{p}$ & \multicolumn{1}{c}{$\bolds{\delta}$} & \multicolumn{1}{c}{$\bolds{\rho}$}
& \multicolumn{1}{c}{$\bolds{S_{n,m}}$}
& \multicolumn{1}{c}{\textbf{RMD test}} & \multicolumn{1}{c}{\textbf{TMD test}} &
\multicolumn{1}{c}{$\bolds{I_{n,m}}$} & \multicolumn{1}{c@{}}{$\bolds{I_{n,m}^{+}}$}\\
\hline
3 & $\bolds{\delta}_{1}$ & -0.25 & 0.79 (90\%) & 0.62 (71\%) &
0.88 &
0.89 (100\%) & 0.89 (100\%)\\
3 & $\bolds{\delta}_{1}$ & 0.00 & 0.64 (82\%) & 0.45 (57\%) & 0.78
& 0.68
(87\%) & 0.68 (87\%)\\
3 & $\bolds{\delta}_{1}$ & 0.25 & 0.53 (78\%) & 0.38 (56\%) & 0.68
& 0.54
(79\%) & 0.54 (79\%)\\
3 & $\bolds{\delta}_{1}$ & 0.50 & 0.51 (73\%) & 0.41 (59\%) & 0.70
& 0.47
(67\%) & 0.47 (67\%)\\
3 & $\bolds{\delta}_{1}$ & 0.90 & 0.62 (64\%) & 0.85 (99\%) & 0.97
& 0.40
(41\%) & 0.41 (42\%)\\
[4pt]
5 & $\bolds{\delta}_{2}$ & -0.25 & 0.93 (95\%) & 0.74 (76\%) &
0.98 &
0.97 (99\%) & 0.97 (99\%)\\
5 & $\bolds{\delta}_{2}$ & 0.00 & 0.80 (87\%) & 0.56 (60\%) & 0.92
& 0.86
(93\%) & 0.86 (93\%)\\
5 & $\bolds{\delta}_{2}$ & 0.25 & 0.59 (73\%) & 0.39 (47\%) & 0.81
& 0.66
(81\%) & 0.66 (81\%)\\
5 & $\bolds{\delta}_{2}$ & 0.50 & 0.56 (67\%) & 0.42 (50\%) & 0.84
& 0.48
(57\%) & 0.48 (57\%)\\
5 & $\bolds{\delta}_{2}$ & 0.90 & 0.63 (64\%) & 0.88 (89\%) & 0.99
& 0.40
(40\%) & 0.40 (40\%)\\
[4pt]
3 & $\bolds{\delta}_{3}$ & -0.25 & 0.74 (89\%) & 0.54 (64\%) &
0.83 &
0.83 (100\%) & 0.83 (100\%)\\
3 & $\bolds{\delta}_{3}$ & 0.00 & 0.56 (87\%) & 0.34 (53\%) & 0.64
& 0.59
(92\%) & 0.59 (92\%)\\
3 & $\bolds{\delta}_{3}$ & 0.25 & 0.42 (87\%) & 0.23 (48\%) & 0.49
& 0.46
(93\%) & 0.46 (93\%)\\
3 & $\bolds{\delta}_{3}$ & 0.50 & 0.33 (86\%) & 0.15 (40\%) & 0.38
& 0.37
(97\%) & 0.37 (97\%)\\
3 & $\bolds{\delta}_{3}$ & 0.90 & 0.27 (83\%) & 0.12 (38\%) & 0.32
& 0.27
(83\%) & 0.27 (83\%)\\
[4pt]
5 & $\bolds{\delta}_{4}$ & -0.25 & 0.92 (95\%) & 0.65 (68\%) &
0.96 &
0.95 (99\%) & 0.95 (99\%)\\
5 & $\bolds{\delta}_{4}$ & 0.00 & 0.75 (90\%) & 0.43 (51\%) & 0.83
& 0.82
(99\%) & 0.82 (99\%)\\
5 & $\bolds{\delta}_{4}$ & 0.25 & 0.49 (87\%) & 0.20 (35\%) & 0.57
& 0.60
(100\%) & 0.60 (100\%)\\
5 & $\bolds{\delta}_{4}$ & 0.50 & 0.41 (90\%) & 0.16 (34\%) & 0.45
& 0.43
(100\%) & 0.43 (100\%)\\
5 & $\bolds{\delta}_{4}$ & 0.90 & 0.29 (92\%) & 0.10 (32\%) & 0.31
& 0.33
(100\%) & 0.33 (100\%)\\
\hline
\end{tabular*}
\end{table}

\setcounter{table}{4}
\begin{table}
\caption{(Continued)}
\begin{tabular*}{\tablewidth}{@{\extracolsep{\fill}}lcd{2.2}llcll@{}}
\hline
\multicolumn{8}{@{}c@{}}{\textbf{Power and RE \%} $\bolds{(n=m=25)}$}\\
\hline
\multicolumn{3}{@{}c}{\textbf{Set up}}
& \multicolumn{3}{c}{\textbf{Directional tests}} &
\multicolumn{2}{c@{}}{\textbf{Integral tests}}\\[-4pt]
\multicolumn{3}{@{}c}{\hrulefill}
& \multicolumn{3}{c}{\hrulefill} &
\multicolumn{2}{c@{}}{\hrulefill}\\
$\bolds{p}$ & \multicolumn{1}{c}{$\bolds{\delta}$} & \multicolumn{1}{c}{$\bolds{\rho}$}
& \multicolumn{1}{c}{$\bolds{S_{n,m}}$}
& \multicolumn{1}{c}{\textbf{RMD test}} & \multicolumn{1}{c}{\textbf{TMD test}} &
\multicolumn{1}{c}{$\bolds{I_{n,m}}$} & \multicolumn{1}{c@{}}{$\bolds{I_{n,m}^{+}}$}\\
\hline
3 & $\bolds{\delta}_{1}$ & -0.25 & 0.96 (98\%) & 0.90 (91\%) &
0.98 &
0.98 (100\%) & 0.98 (100\%)\\
3 & $\bolds{\delta}_{1}$ & 0.00 & 0.85 (92\%) & 0.72 (78\%) & 0.92
& 0.85
(92\%) & 0.86 (92\%)\\
3 & $\bolds{\delta}_{1}$ & 0.25 & 0.80 (88\%) & 0.69 (76\%) & 0.90
& 0.75
(83\%) & 0.75 (83\%)\\
3 & $\bolds{\delta}_{1}$ & 0.50 & 0.75 (84\%) & 0.67 (75\%) & 0.89
& 0.66
(74\%) & 0.66 (74\%)\\
3 & $\bolds{\delta}_{1}$ & 0.90 & 0.89 (89\%) & 0.98 (99\%) & 1.00
& 0.59
(59\%) & 0.61 (61\%)\\
[4pt]
5 & $\bolds{\delta}_{2}$ & -0.25 & 1.00 (100\%) & 0.98 (98\%) &
1.00 &
1.00 (100\%) & 1.00 (100\%)\\
5 & $\bolds{\delta}_{2}$ & 0.00 & 0.96 (97\%) & 0.85 (86\%) & 0.99
& 0.98
(99\%) & 0.98 (99\%)\\
5 & $\bolds{\delta}_{2}$ & 0.25 & 0.85 (88\%) & 0.74 (76\%) & 0.97
& 0.83
(86\%) & 0.83 (86\%)\\
5 & $\bolds{\delta}_{2}$ & 0.50 & 0.81 (84\%) & 0.74 (77\%) & 0.96
& 0.70
(73\%) & 0.70 (73\%)\\
5 & $\bolds{\delta}_{2}$ & 0.90 & 0.90 (90\%) &0.99
(100\%) &
1.00 & 0.57 (57\%) & 0.58 (58\%)\\
[4pt]
3 & $\bolds{\delta}_{3}$ & -0.25 & 0.94 (96\%) & 0.85 (87\%) &
0.98 &
0.96 (98\%) & 0.96 (98\%)\\
3 & $\bolds{\delta}_{3}$ & 0.00 & 0.75 (92\%) & 0.57 (69\%) & 0.82
& 0.79
(96\%) & 0.79 (96\%)\\
3 & $\bolds{\delta}_{3}$ & 0.25 & 0.62 (89\%) & 0.39 (56\%) & 0.70
& 0.66
(94\%) & 0.66 (94\%)\\
3 & $\bolds{\delta}_{3}$ & 0.50 & 0.54 (90\%) & 0.31 (52\%) & 0.60
& 0.55
(92\%) & 0.55 (92\%)\\
3 & $\bolds{\delta}_{3}$ & 0.90 & 0.44 (90\%) & 0.20 (42\%) & 0.49
& 0.42
(86\%) & 0.42 (86\%)\\
[4pt]
5 & $\bolds{\delta}_{4}$ & -0.25 & 0.99 (99\%) & 0.94 (94\%) &
1.00 &
1.00 (100\%) & 1.00 (100\%)\\
5 & $\bolds{\delta}_{4}$ & 0.00 & 0.94 (96\%) & 0.72 (74\%) & 0.97
& 0.97
(100\%) & 0.97 (100\%)\\
5 & $\bolds{\delta}_{4}$ & 0.25 & 0.71 (90\%) & 0.41 (52\%) & 0.79
& 0.79
(100\%) & 0.79 (100\%)\\
5 & $\bolds{\delta}_{4}$ & 0.50 & 0.58 (91\%) & 0.25 (39\%) & 0.63
& 0.64
(100\%) & 0.64 (100\%)\\
5 & $\bolds{\delta}_{4}$ & 0.90 & 0.42 (87\%) & 0.18 (36\%) & 0.49
& 0.46
(94\%) & 0.46 (94\%)\\
\hline
\end{tabular*}
\end{table}

As expected, in every case the power of the TMD test is higher than
that of
$S_{n,m}$ test and the RMD test. The $S_{n,m}$ test is almost always more
powerful than the RMD test. The relative efficiency of $S_{n,m}$
compared to
the TMD test is quite high in most cases. When $n=m=15$ the relative
efficiency ranges between 65--95\%. It is almost always above 90\% when the
sample size increases to 25 per group. In general the two integral
tests had
very similar power. They had larger power than $S_{n,m}$ when $\rho<0$. As
$\rho$ increased, the power of $S_{n,m}$ improved relative to the two integral
tests. Test (\ref{Snm}) seems to perform better when the components of
$\bolds{\delta}$ were unequal. We also note that when the integral tests
outperform $S_{n,m}$ the difference is usually small, whereas the $S_{n,m}$
test can outperform the integral tests substantially. For example, observe
pattern 2 where the powers of $S_{n,m}$ and $I_{n,m}$ are $0.93$ and
$0.97$,
respectively, when $\rho=-0.25$ and $0.63$ versus $0.40 $ when
$\rho=0.90$.

\section{Illustration}\label{sec5}

Prior to conducting a two-year rodent cancer bioassay to evaluate the
toxicity/carcinogenicity of a chemical, the National Toxicology Program (NTP)
routinely conducts a 90-day pre-chronic dose finding study. One of the goals
of the 90-day study is to determine the maximum tolerated dose (MTD)
that can
be used in the two-year chronic exposure study. Accurate determination
of the
MTD is critical for the success of the two-year cancer bioassay. Cancer
bioassays are typically very expensive and time consuming. Therefore their
proper design, that is, choosing the correct dosing levels, is very important.
When the highest dose used in the two-year study exceeds the MTD, a large
proportion of animals in the high dose group(s) may die well before the
end of
the study, and the data from such group(s) cannot be used reliably. This
results in inefficiency and wasted resources.

Typically the NTP uses the 90-day study to determine the MTD on the
basis of
a large number of correlated endpoints that provide information regarding
toxicity. These include body weight, organ weights, clinical chemistry (red
blood cell counts, cell volume, hemoglobin, hematocrit, lymphocytes, etc.),
histopathology (lesions in various target organs), number of deaths and so
forth. The dose response data is analyzed for each variable separately using
Dunnett's or the Williams's test (or their nonparametric versions,
Dunn's test
and Shirley's test, resp.). NTP combines results from all such analyses
qualitatively and uses other biological and toxicological information when
making decisions regarding the highest dose for the two-year cancer bioassay.
Analyzing correlated variables one at a time may result in loss of
information. The proposed methodology provides a convenient method to combine
information from several outcome variables to make comparisons between groups.

We now illustrate our methodology by re-analyzing data obtained from a recent
NTP study of the chemical Citral [\citet{autokey28}]. Citral is a flavoring
agent that
is widely used in a variety of food items. The NTP assigned a random
sample of
10 male rats to the control group and 10 to the 1785 mg/kg dose
group. Hematological and clinical chemistry measurements such as the
number of
platelets (in 1000~per l), urea nitrogen (UN) (in mg/dl),
alkaline phosphatase (AP) (in IU/l) and bile acids (BA) (in
mol/l) were recorded on each animal at the end of the study.
The NTP
performed univariate analysis on each of these variables and found no
significant difference between the control and dose group except for the
concentration of urea nitrogen which was increased in the high dose group.
This increase was marginally significant at the $5\%$ level and not at all
after correcting for multiplicity. We applied the proposed methodology to
compare the control with the high-dose group (1785 mg/kg) in terms of all
nonnegative linear combinations of the above mentioned four variables. We
test the null hypothesis of no difference between the control and the
high-dose group against the alternative that the high-dose group is
stochastically
larger (in the above four variables) than the control group. The resulting
$p$-value based on 10,000 bootstrap samples was $0.025$, which is significant
at a $5\%$ level of significance. The estimated value of $\mathbf{s}
_{\max}$ was $(0.074,0.986,0.012,0.150)^{T}$ and the estimated 95\% confidence
region is given by $\{\mathbf{s}\in\mathcal{S}_{+}^{p-1}\dvtx \widehat
{\mathbf{s}}^{T}_{\max}\mathbf{s}\leq0.93\}$. Hence the
confidence set
includes any $\mathbf{s}$ which is within $21.5^{\circ}$ degrees of
$\widehat{\mathbf{s}}_{\max}$. This is a relatively large set due
to the
small sample sizes. Clearly our methodology appears to be sensitive to detect
statistical differences which were not noted by NTP. Furthermore, our
methodology allows us to infer that indeed 1785 mg/kg dose group is
larger in
the multivariate stochastic order than the control group. This is a much
stronger conclusion than the simple ordering of their means. Thus we believe
that the proposed framework and methodology for studying ordered distributions
can serve as a useful tool in toxicology and is also applicable to a wide
range of other problems as alluded to in this paper.

\section{Concluding remarks and some open problems}\label{sec6}

In many applications, researchers are interested in comparing two experimental
conditions, for example, a treatment and a control group, in terms of a
multivariate
response. In classical multivariate analysis one addresses such
problems by
comparing the mean vectors using Hotelling's $T^{2}$ statistic. The assumption
of MVN, underlying Hotelling's $T^{2}$ test, may not hold in practice.
Moreover if the data is not MVN, then the comparison of population
means may
not always provide complete information regarding the differences
between the
two experimental groups. Secondly, Hotelling's $T^{2} $ statistics are designed
for two-sided alternatives and may not be ideal if a researcher is interested
in one-sided, that is, ordered alternatives. Addressing such problems requires
one to compare the two experimental groups nonparametrically in terms
of the
multivariate stochastic order. Such comparisons, however, are very high
dimensional and not easy to perform.

In this article we circumvent this challenge by considering the notion
of the
linear stochastic order between two random vectors. The linear stochastic
order is a ``weak'' generalization of the univariate stochastic order. The
linear stochastic order is simple to interpret and has an intuitive appeal.
Using this notion of ordering, we developed nonparametric directional
inference procedures. Intuitively, the proposed methodology seeks to determine
the direction that best separates two multivariate populations. Asymptotic
properties of the estimated direction are derived. Our test based on
the best
separating direction may be viewed as a generalization of Roy's classical
largest root test for comparing several MVN populations. To the best
of our
knowledge this is the first general test for multivariate ordered
distributions. Since in practice sample sizes are small, we use the bootstrap
methodology for drawing inferences.

We illustrated the proposed methodology using a data obtained from a recent
toxicity/carcinogenicity study conducted by the US National Toxicology Program
(NTP) on the chemical Citral. A re-analysis of their 90-day data using our
proposed methodology revealed a linear stochastic increase in
platelets, urea
nitrogen, alkaline phosphatase and bile acids in the high-dose group relative
to the control group, which was not seen in the original univariate analysis
conducted by the NTP. These findings suggest that the proposed
methodology may
have greater sensitivity than the commonly used univariate statistical
procedures. Our methodology is sufficiently general since it is nonparametric
and can be applied to discrete and/or continuous outcome variables.
Furthermore, our methodology exploits the underlying dependence
structure in
the data, rather than analyzing one variable at a time.

We note that our example and some of our results pertain to continuous RVs.
However, the methodology may be used, with appropriate modification (e.g.,
methods for dealing with ties) with discrete (or mixed) data with no problem.
Although the focus of this paper has been the comparison of two multivariate
vectors, in many applications, especially in dose response studies,
researchers may be interested in determining trends (order) among several
groups. Similarly to classical parametric order restricted inference literature,
one could generalize the methodology developed in this paper to test
for order
restrictions among multiple populations. For example, one could extend the
results to $K\geq2$ RVs ordered by the simple ordering, that is,
$\mathbf{X}%
_{1}\prec_{\lst }\mathbf{X}_{2}\prec_{\lst }\cdots\prec_{\lst }%
\mathbf{X}_{K}$ or to RVs ordered by the tree ordering, that is,
$\mathbf{X}_{1}\prec_{\lst }\mathbf{X}_{j}$ where $j=2,\ldots,K$. As
pointed out by a referee the hypotheses $H_{0}\dvtx \mathbf{X}\preceq
_{\st }\mathbf{Y}$ versus $H_{1}\dvtx \mathbf{X}\nprec_{\st }\mathbf{Y}$
can also be formulated and tested using the approach described. First note
that the null hypothesis implies $\Psi(\mathbf{s})\geq1/2$ for all
$\mathbf{s}\in\mathcal{S}_{+}^{p-1}$. On the other hand under the
alternative there is an $\mathbf{s}\in\mathcal{S}_{+}^{p-1}$ for which
$\Psi(\mathbf{s})<1/2$. Thus a test may be based on the statistic%
\[
N^{1/2}\bigl(\Psi_{n,m}(\widehat{\mathbf{s}}_{\min})-1/2
\bigr),
\]
where $\widehat{\mathbf{s}}_{\min}$ is the value which minimizes
$\Psi_{n,m}(\mathbf{s})$. It is also clear that the least favorable
configuration occurs when $\Psi(\mathbf{s})=1/2$ for all
$\mathbf{s}%
\in\mathcal{S}_{+}^{p-1}$ which is equivalent to $\mathbf{X}%
=_{\st }\mathbf{Y}$.

We believe that the result obtained here may be useful beyond order restricted
inference. Our simulation study suggests that our estimator of the best
separating direction, that is, (\ref{s-max-hat}) may be useful even in the
context of classical multivariate analysis where it may be viewed as a robust
alternative to Roy's classical estimate. Finally we note that the linear
stochastic order may be useful in a variety of other statistical
problems. For
example, we believe that it provides a useful framework for linearly combining
the results of several diagnostic markers. This is a well-known problem
in the
context of ROC curve analysis in diagnostic medicine.

\begin{appendix}\label{app}
\section*{Appendix: Proofs}

\begin{pf*}{Proof of Theorem~\ref{Thm-Closure}}
(i) Let $g\dvtx
\mathbb{R}^{p}\rightarrow
\mathbb{R}^{n}$ be an affine increasing function. Clearly $g ( \mathbf{x}
) =\mathbf{v}+\mathbf{M}\mathbf{x}$ for some $n$ vector
$\mathbf{v}$ and $n\times p$ matrix $\mathbf{M}$ with nonnegative
elements. Thus for any $\mathbf{u}\in$ $
\mathbb{R}_{+}^{n}$ we have $\mathbf{s}=\mathbf{M}^{T}\mathbf{u}\in
$ $
\mathbb{R}_{+}^{p}$. Hence
\[
\mathbf{u}^{T}g(\mathbf{X})\!=\!\mathbf{u}^{T}(\mathbf{v}%
\!+\!\mathbf{MX})\!=\!\mathbf{u}^{T}\mathbf{v}\!+\!\mathbf{s}^{T}
\mathbf{X}\preceq_{\st }\mathbf{u}^{T}\mathbf{v}\!+\!
\mathbf{s}^{T}\mathbf{Y}\!=\!\mathbf{u}^{T}(\mathbf{v}\!+\!\mathbf{MY})
\!=\!\mathbf{u}^{T}g(\mathbf{Y})
\]
as required where the inequality holds because $\mathbf{X}\preceq
_{\lst }\mathbf{Y}$. (ii) Fix $I\in\{ 1,\ldots,\break p
\}$. Let $\mathbf{X}=(\mathbf{X}_{I},\mathbf{X}%
_{\overline{I}})$, $\mathbf{Y}=(\mathbf{Y}_{I},\mathbf{Y}%
_{\overline{I}})$ where $\overline{I}$ is the complement of $I$ in
$ \{
1,\ldots,p \}$. Further define $\mathbf{s}^{T}=(\mathbf{s}
_{I}^{T},\mathbf{s}_{\overline{I}}^{T})$ where $\mathbf{s}\in$ $
\mathbb{R}_{+}^{p}$, and set $\mathbf{s}_{\overline{I}}^{T}=0$. It follows
that for
all $\mathbf{s}_{I}\in
\mathbb{R}^{\dim( I ) }$ we have%
\[
\mathbf{s}_{I}^{T}\mathbf{X}_{I}=
\mathbf{s}^{T}\mathbf{X}%
\preceq_{\st }\bolds{
\mathbf{s}^{T}\mathbf{Y}=s}_{I}^{T}%
\mathbf{Y}_{I}%
\]
as required. (iii) Let $\phi\dvtx
\mathbb{R}
\rightarrow
\mathbb{R}
$ be any increasing function. Note that%
\[
\mathbb{E} \bigl( \phi\bigl( \mathbf{s}^{T}\mathbf{X} \bigr) \bigr) =
\mathbb{E\bigl(E} \bigl( \phi\bigl( \mathbf{s}^{T}\mathbf{X}%
\bigr) |\mathbf{Z} \bigr) \bigr)\leq\mathbb{E\bigl(E} \bigl( \phi\bigl(
\mathbf{s}^{T}\mathbf{Y} \bigr) |\mathbf{Z} \bigr) \bigr)=\mathbb{E} \bigl(
\phi\bigl( \mathbf{s}^{T}\mathbf{Y} \bigr) \bigr).
\]
The inequality is a consequence of $\mathbf{X}|\mathbf{Z}%
=\mathbf{z}\preceq_{\lst }\mathbf{Y}|\mathbf{Z}=\mathbf{z}$.
Since $\phi$ is arbitrary it follows that $\mathbf{X}\preceq_{\lst }\mathbf
{Y}$ as required. (iv) Let $\mathbf
{X}%
=(\mathbf{X}_{1},\ldots,\mathbf{X}_{n})$, and define
$\mathbf{Y}$
similarly. Let $\mathbf{s}\in$ $
\mathbb{R}_{+}^{p}$ where $p=p_{1}+\cdots+p_{n}$. Now%
\[
\mathbf{s}^{T}\mathbf{X}=\mathbf{s}_{1}^{T}
\mathbf{X}_{1}+\cdots+\mathbf{s}_{n}^{T}
\mathbf{X}_{n}\quad\mbox{and}\quad\mathbf{s}%
^{T}
\mathbf{Y}=\mathbf{s}_{1}^{T}\mathbf{Y}_{1}+\cdots+
\mathbf{s}_{n} 
^{T}\mathbf{Y}_{n}%
\]
by assumption $\mathbf{s}_{i}^{T}\mathbf{X}_{i}\preceq_{\st }%
\mathbf{s}_{i}^{T}\mathbf{Y}_{i}$ for $i=1,\ldots,n$. In addition
$\mathbf{s}_{i}^{T}\mathbf{X}_{i}$ and $\mathbf{s}_{j}%
^{T}\mathbf{X}_{j}$ are independent for $i\neq j$. It follows from Theorem
1.A.3 in \citet{ShaSha07} that $\mathbf{s}_{1}%
^{T}\mathbf{X}_{1}+\cdots+\mathbf{s}_{n}^{T}\mathbf{X}_{n}%
\preceq_{\st }\mathbf{s}_{1}^{T}\mathbf{Y}_{1}+\cdots+\mathbf
{s}%
_{n}^{T}\mathbf{Y}_{n}$, that is, $\mathbf{X}\preceq_{\lst }%
\mathbf{Y}$ as required. (v) By assumption $\mathbf{X}%
_{n}\Rightarrow\mathbf{X}$ and $\mathbf{Y}_{n}\Rightarrow
\mathbf{Y}$ where the symbol $\Rightarrow$ denotes convergence in
distribution. By the continuous mapping theorem $\mathbf{s}^{T}%
\mathbf{X}_{n}\Rightarrow\mathbf{s}^{T}\mathbf{X}$ and
$\mathbf{s}^{T}\mathbf{Y}_{n}\Rightarrow\mathbf{s}^{T}%
\mathbf{Y}$. It follows that
%
\begin{equation}\label{Pf-CL-1}
\mathbb{P} \bigl( \mathbf{s}^{T}\mathbf{X}_{n}\geq t \bigr)
\rightarrow\mathbb{P} \bigl( \mathbf{s}^{T}\mathbf{X}\geq t \bigr)
\quad\mbox{and}\quad\mathbb{P} \bigl( \mathbf{s}^{T}\mathbf{Y}_{n}\geq t
\bigr) \rightarrow\mathbb{P} \bigl( \mathbf{s}^{T}\mathbf{Y}\geq t
\bigr).%
\end{equation}
Moreover since $\mathbf{X}_{n}\preceq_{\lst }\mathbf{Y}_{n}$ we
have%
%
\begin{equation}
\label{Pf-CL-2} \mathbb{P} \bigl( \mathbf{s}^{T}\mathbf{X}_{n}
\geq t \bigr) \leq\mathbb{P} \bigl( \mathbf{s}^{T}\mathbf{Y}_{n}
\geq t \bigr) \qquad\mbox{for all }n\in\mathbb{N}.
\end{equation}
Combining (\ref{Pf-CL-1}) and (\ref{Pf-CL-2}) we have $\mathbb{P} (
\mathbf{s}^{T}\mathbf{X}\geq t ) \leq\mathbb{P} (
\mathbf{s}^{T}\mathbf{Y}\geq t ) $, that is, $\mathbf
{X}%
\preceq_{\lst }\mathbf{Y}$ as required.
\end{pf*}

Before proving Theorem~\ref{Thm-ERVs}, we provide a definition and a
preliminary lemma.
%
\begin{definition}
We say that the RV $\mathbf{X}$ has an elliptical distribution with
parameters $\bolds{\mu}$ and $\bolds{\Sigma}$ and generator
$\phi( \cdot) $, denoted $\mathbf{X}\sim E_{p} (
\bolds{\mu},\bolds{\Sigma},\phi)$, if its characteristic
function is given by $\exp( i\mathbf{t}^{T}\bolds{\mu
} )
\phi( \mathbf{t}^{T}\bolds{\Sigma} \mathbf{t} )$.
\end{definition}

For this and other facts about elliptical distributions which we use in the
proofs below, see \citet{FanKotNg89}.
%
\begin{lemma}
\label{Lemma-Elliptic}Let $X\sim E_{1} ( \mu,\sigma,\phi) $
and $Y\sim E_{1} ( \mu^{\prime},\sigma^{\prime},\phi) $ be
univariate elliptical RVs supported on $ \mathbb{R}$. Then
$X\preceq_{\st }Y$ if and only if $\mu\leq\mu^{\prime}$ and
$\sigma=\sigma^{\prime}$.
\end{lemma}
\begin{pf}
Since $X$ and $Y$ have the same generator they have the stochastic
representation:%
%
\begin{equation}\label{a}
X=_{\st }\mu+\sigma RU\quad\mbox{and}\quad Y=_{\st }\mu^{\prime}+
\sigma^{\prime}RU,
\end{equation}
where $R$ is a nonnegative RV, independent of the RV $U$, satisfying
$\mathbb{P}(U=\pm1)=1/2$; cf. \citet{FanKotNg89}. It follows that
$RU$ is a symmetric RV supported on $\mathbb{R} $ with a strictly
increasing DF which we denoted by $F_{0}$. Let $F_{X}$ and $F_{Y}$
denote the DFs of $X$ and $Y$, respectively. Note that $X\preceq_{\st }Y$
if and only if $F_{X} ( t ) \geq$ $F_{Y} ( t ) $ for
all $t\in\mathbb{R}
$, or equivalently by\vadjust{\goodbreak} (\ref{a}), if and only if%
%
\begin{equation}\label{1}
F_{0}\biggl(\frac{t-\mu}{\sigma}\biggr)\geq F_{0}\biggl(
\frac{t-\mu^{\prime}}{\sigma^{\prime}}\biggr)
\end{equation}
for all $t\in\mathbb{R}$. It is obvious that (\ref{1}) holds when
$\mu\leq\mu^{\prime}$ and $\sigma=\sigma^{\prime}$, establishing
sufficiency. Now assume that $X\preceq_{\st }Y$. Put $t=\mu$ in
(\ref{1}), and use the strict monotonicity of $F_{0}$ to get
$0\geq(\mu-\mu^{\prime})/\sigma^{\prime}$, that is, $\mu^{\prime
}\geq\mu$. Suppose now that $\sigma^{\prime}>\sigma$. It follows from
(\ref{1}) and the the strict monotonicity of $F_{0}$ that $(t-\mu
)/\sigma\geq(t-\mu^{\prime})/\sigma^{\prime}$ which is equivalent to
$t\geq(
\mu\sigma^{\prime}-\mu^{\prime}\sigma) / ( \sigma^{\prime}%
-\sigma)$. The latter, however, contradicts the fact that
(\ref{1}) holds for all $t\in\mathbb{R}$. A similar argument shows
that $\sigma^{\prime}<\sigma$ cannot hold; hence we must have
$\sigma=\sigma^{\prime}$ as required.
\end{pf}
%
\begin{remark}
Note that Lemma~\ref{Lemma-Elliptic} may not hold for distributions
with a finite support. For example, if $R\sim U ( 0,1 ), $ then
by (\ref{a}) $X$ $\sim U ( \mu-\sigma,\mu+\sigma) $ and
$Y\sim U (
\mu^{\prime}-\sigma^{\prime},\mu^{\prime}+\sigma^{\prime} )$. It
is easily verified that in this case $X\preceq_{\st }Y$ if and only if
$\Delta=\mu^{\prime}-\mu\geq0$ and $-\Delta\leq\sigma^{\prime}-\sigma
\leq\Delta$; that is, it is not required that $\sigma=\sigma^{\prime}$.
Hence the assumption that $X$ and $Y$ are supported on $\mathbb{R} $ is
necessary.
\end{remark}

We continue with the proof of Theorem~\ref{Thm-ERVs}.
\begin{pf*}{Proof of Theorem~\ref{Thm-ERVs}}
Let $\mathbf{X}$ and $\mathbf{Y}$ be be $E_{p} (
\bolds{\mu},\bolds{\Sigma},\phi) $ and $E_{p} (
\bolds{\mu}^{\prime},\break \bolds{\Sigma}^{\prime},\phi) $
supported on $ \mathbb{R}^{p}$. Suppose that
$\mathbf{X}\preceq_{\lst }\mathbf{Y}$. Choose
$\mathbf{s}=\mathbf{e}_{i}$ where $\mathbf{e}_{ik}=1$ if $i=k$ and $0$
otherwise. It now follows from Definition~\ref{Def-lst} that $X_{i}%
\preceq_{\st }Y_{i}$. Since $X_{i}$ and $Y_{i}$ are marginally
elliptically distributed RVs with the same generator and supported on~$
\mathbb{R}, $ then by Lemma~\ref{Lemma-Elliptic} we must have
%
\begin{equation}\label{m-1}
\mu_{i}\leq\mu_{i}^{\prime}\quad\mbox{and}\quad
\sigma_{ii}=\sigma_{ii}^{\prime}.
\end{equation}
The latter holds, of course, for all $1\leq i\leq p$. Choosing
$\mathbf{s}=\mathbf{e}_{i}+\mathbf{e}_{j}$ we have $X_{i}+X_{j}\preceq
_{\st }Y_{i}+Y_{j}$. Note that $X_{i}+X_{j}$ and $Y_{i}+Y_{j}$ are
supported on $ \mathbb{R} $ and follow a univariate elliptical
distribution with the same generator [\citet{FanKotNg89}]. Applying
Lemma~\ref{Lemma-Elliptic} again we find that
%
\begin{equation}\label{m-2}
\mu_{i}+\mu_{j}\leq\mu_{i}^{\prime}+
\mu_{j}^{\prime}\quad\mbox{and}\quad\sigma_{ii}+
\sigma_{jj}+2\sigma_{ij}=\sigma_{ii}^{\prime}+
\sigma_{jj}^{\prime
}+2\sigma_{ij}^{\prime}.
\end{equation}
The latter holds, of course, for all $1\leq i\neq j\leq p$. It is easy
to see
that equations (\ref{m-1}) and (\ref{m-2}) imply that $\bolds{\mu}%
\leq\bolds{\mu}^{\prime}$ and $\bolds{\Sigma=\Sigma}^{\prime}$.
Recall [cf. \citet{FanKotNg89}] that we may write $\mathbf{X}=_{\st }%
\bolds{\mu}+R\mathbf{SU}$ and $\mathbf{Y}=_{\st }\bolds
{\mu}^{\prime}+R\mathbf{SU}$ where
$\bolds{\Sigma=S}^{T}\mathbf{S,}$ $\mathbf{U}$~is a
uniform RV on $\mathcal{S}_{+}^{p-1}$, and $R$ is a nonnegative RV.
Let $S$ be an upper set in~$\mathbb{R}^{p}$. Clearly the set $[S-\bolds
{\mu}]:= \{ \mathbf{x}%
-\bolds{\mu}\dvtx x\in S \} $ is also an upper set and
$[S-\bolds{\mu}]\subseteq[ S-\bolds{\mu}^{\prime}]$ since
$\bolds{\mu}\leq\bolds{\mu}^{\prime}$. Now,
\[
\mathbb{P} ( \mathbf{X}\in S ) =\mathbb{P} \bigl( \mathbf{X}_{0}\in
[ S-\bolds{\mu}] \bigr) \leq\mathbb{P} \bigl( \mathbf{X}_{0}\in
\bigl[ S-\bolds{\mu}^{\prime}\bigr] \bigr) =\mathbb{P} ( \mathbf{Y}\in S ),
\]
where $\mathbf{X}_{0}=R\mathbf{SU}$, hence $\mathbf
{X}\preceq_{\st }\mathbf{Y}$. This proves the ``if'' part. The ``only if'' part follows
immediately.\vadjust{\goodbreak}
\end{pf*}
\begin{pf*}{Proof of Theorem~\ref{Thm-MVBs}}
Let $\mathcal{X}_{p}=\{\mathbf{x}\dvtx  ( x_{1},\ldots,x_{p} )
\in\{ 0,1 \}^{p}\}$ denote the support of a $p$-dimensional
multivariate binary (MVB) RV. By definition the relationship
$\mathbf{X}%
\preceq_{\lst }\mathbf{Y}$ implies that for all $(t,\mathbf
{s})\in\mathbb{R}_{+}\times\mathbb{R}_{+}^{p}$,
%
\begin{equation}\label{mvb1}
\mathbb{P} \bigl( \mathbf{s}^{T}\mathbf{X}>t \bigr) \leq\mathbb{P}
\bigl( \mathbf{s}^{T}\mathbf{Y}>t \bigr).%
\end{equation}
Now note that%
%
\begin{equation}\label{mvb2}\hspace*{28pt}
\mathbb{P} \bigl( \mathbf{s}^{T}\mathbf{X}>t \bigr) =\sum
_{\mathbf{x}\in\mathcal{X}_{p}}f ( \mathbf{x} ) \mathbb{I}_{(\mathbf
{s}^{T}\mathbf{x}>t)}\quad\mbox{and}\quad
\mathbb{P} \bigl( \mathbf{s}^{T}\mathbf{Y}>t \bigr) =\sum
_{\mathbf{x}\in
\mathcal{X}_{p}}g ( \mathbf{x} ) \mathbb{I}_{(\mathbf
{s}%
^{T}\mathbf{x}>t)},
\end{equation}
where $f$ and $g$ are the probability mass functions of $\mathbf
{X}$ and
$\mathbf{Y}$, respectively. Let $U$ be an upper set on $\mathcal{X}_{p}$.
It is well known [cf. \citet{DavPri02}] that $U$ can be
written as
%
\begin{equation}\label{mvb3}
U=\bigcup_{j\in J}U ( \mathbf{x}_{j} ),
\end{equation}
where $\mathbf{x}_{j}$ are the distinct minimal elements of $U$, and
$U ( \mathbf{x}_{j} ) = \{ \mathbf{x}\dvtx \mathbf
{x}%
\geq\mathbf{x}_{j} \} $ are themselves upper sets [in fact
$U (
\mathbf{x}_{j} ) $ is an upper orthant]. The set $\{\mathbf
{x}%
_{j}\dvtx $ $j\in J\}$ is often referred to as an anti-chain. Now observe
that for any $\mathbf{s}\in\mathbb{R}_{+}^{p}$ the set $ \{
\mathbf{x}\dvtx \mathbf
{s}^{T}\mathbf{x}%
>t \} $ is an upper set. Hence it must be of the form of (\ref
{mvb3}) for some anti-chain $\{\mathbf{x}_{j}\dvtx $ $j\in J\}$.
Suppose now, that for some $U\in\mathcal{X}_{p}$ there is a vector
$\mathbf{s}_{U}\in\mathbb{R}_{+}^{p}$ such that $U= \{ \mathbf{x}\dvtx \mathbf
{s}_{U}%
^{T}\mathbf{x}>t \} $ for some fixed $t>0$. Then using (\ref{mvb1})
and (\ref{mvb2}) we have%
\begin{eqnarray*}
\mathbb{P} ( \mathbf{X}\in U ) &=& \sum_{\mathbf{x}%
\in\{ \mathbf{x}\dvtx \mathbf{s}_{U}^{T}\mathbf{x}>t \}
}f (
\mathbf{x} ) =\mathbb{P} \bigl( \mathbf{s}_{U}%
^{T}
\mathbf{X}>t \bigr) \leq\mathbb{P} \bigl( \mathbf{s}%
_{U}^{T}
\mathbf{Y}>t \bigr)\\
&=&\sum_{\mathbf{x}\in\{
\mathbf{x}\dvtx \mathbf{s}_{U}^{T}\mathbf{x}>t \} }g ( \mathbf{x} ) =
\mathbb{P} ( \mathbf{Y}\in U ).
\end{eqnarray*}
We will complete the proof by showing that for each upper set $U\in
\mathcal{X}_{p}$, we can find a vector $\mathbf{s}_{U}$ for which
$\mathbf{s}_{U}^{T}\mathbf{x}>t$ for $\mathbf{x}\in U$ and
$\mathbf{s}_{U}^{T}\mathbf{x}\leq t$ for $\mathbf{x}\in
U^{c}=\mathcal{X}\setminus U$ if and only if $p\leq3$. To do so we
will first
solve the system of equations $\mathbf{s}^{T}\mathbf{x}_{j}=t$ for
$j\in J$. This system can also be written as $\mathbf{X}\mathbf{s}%
=\mathbf{t}$ where
\[
\mathbf{X}=\pmatrix{ \mathbf{x}_{1}
\cr
\vdots
\cr
\mathbf{x}_{J}}
\]
is a $J\times p$ matrix whose rows are the member of the anti-chain defining
$U$, and $\mathbf{t}= ( t,\ldots,t ) $ has dimension $J$.
Clearly the elements of $\mathbf{X}$ are ones and zeros. If $J\leq p$, the
matrix $\mathbf{X}$ is of full rank since its rows are linearly
independent by
the fact that they are an anti-chain. Hence a solution for $\mathbf{s}$
exists. With a bit of algebra, we can further show that a solution
$\mathbf{s}\geq0$ exists. This, of course, is trivially verified when
$p\leq3$. Now set $\mathbf{s}_{U}=\mathbf{s}+\bolds
{\varepsilon}$
for some $\bolds{\varepsilon}\geq0$. It is clear that we can choose
$\bolds{\varepsilon}$ small enough to guarantee that $\mathbf
{s}%
_{U}^{T}\mathbf{x}>t$ if and only if $\mathbf{x}\in U$. Hence if $J\leq
p$, upper set (\ref{mvb3}) can be mapped to a vector $\mathbf{s}_{U}$.
However, the inequality $J\leq p$ for \textit{all} upper sets
$U\subset\mathcal{X}_{p}$ holds if and only if $p\leq3$. This can be
easily shown by enumerating all $18$ upper sets belonging to
$\mathcal{X}_{3}$ [cf. \citet{DavPed11}] and noting that they have
at most three minimal elements. Hence if $p\leq3$, then
$\mathbf{X}\preceq_{\lst }\mathbf{Y}$
$\Longleftrightarrow\mathbf{X}\preceq_{\st }\mathbf{Y}$ as required.

Now let $p=4$, and consider the upper set $U$ generated by the anti-chain
$\mathbf{x}_{j}$, $j=1,\ldots,J$ where $\mathbf{x}_{j}$ are all the
distinct permutations of the vector $ ( 1,1,0,0 )$. Clearly $J=6$.
Note that although $J>p$, the system of equations $\mathbf{X}\mathbf
{s}%
=\mathbf{t}$ is uniquely solved by $\mathbf{s}_{\ast}^{T}= (
t/2,t/2,t/2,t/2 )$. However, this solution coincides with the solution
of the system $\mathbf{X}^{\prime}\mathbf{s}=\mathbf{t}$ where
$\mathbf{X}^{\prime}$ is any matrix obtained from $\mathbf{X}$ by
deleting any
two (or just one) of its rows. Note that the rows of $\mathbf{X}^{\prime}$
correspond to an upper set $U^{\prime}\subset U$. This, in turn,
implies that
for any such $U^{\prime}$ one cannot find a vector $\mathbf{s}%
_{U^{\prime}}$ satisfying $\mathbf{s}_{U^{\prime}}^{T}\mathbf
{x}>t$ if
and only if $\mathbf{x}\in U^{\prime}$ because the inequality will hold
for all $\mathbf{x}\in U$. Thus $U^{\prime}$ does not define an
upper half
plane. This shows that the linear stochastic order and the multivariate
stochastic order do not coincide when $p=4$. A similar argument may be used
for any $p\geq5$. This completes the proof.
\end{pf*}

We first define the term copula.
%
\begin{definition}
Let $F$ be the DF of a $p$-dimensional RV with marginal DFs
$F_{1},\ldots,F_{p}$. The copula $C$ associated with $F$ is a DF such
that%
\[
F ( \mathbf{x} ) =C(\mathbf{x})=C \bigl( F_{1} ( x_{1} ),\ldots,F_{p} ( x_{p} ) \bigr).
\]
It follows that the tail-copula $\overline{C}(\bolds{\cdot})$ is nothing
but the tail of the DF $C(\bolds{\cdot})$.
\end{definition}
\begin{pf*}{Proof of Theorem~\ref{Thm-Copula}}
Suppose that $\mathbf{X}$ and $\mathbf{Y}$ have the same
copula. Let
$\mathbf{X}\preceq_{\lst }\mathbf{Y}$. Choosing $\mathbf
{s}=\mathbf{e}_{i} $
where $\mathbf{e}_{ik}=1$ if $i=k$ and $0$ otherwise, we find using the
definition that $X_{i}\preceq_{\st }Y_{i}$. The latter holds, of course,
for all
$1\leq i\leq p$. Applying Theorem 6.B.14 in \citet{ShaSha07}, we find that $\mathbf{X}\preceq_{\st }\mathbf{Y}$.
The reverse direction is immediate.
\end{pf*}
\begin{pf*}{Proof of Theorem~\ref{Thm-ORT}}
Note that for any $\mathbf{x}\in\mathbb{R}^{p}$ we have%
\begin{eqnarray*}
F ( \mathbf{x} ) &=& C_{\mathbf{X}} \bigl( F_{1}%
(x_{1}),\ldots,F_{p}(x_{p}) \bigr) \geq
C_{\mathbf{X}} \bigl( G_{1}%
(x_{1}),\ldots,G_{p}(x_{p}) \bigr)
\\
&\geq& C_{\mathbf{Y}} \bigl( G_{1}(x_{1}),\ldots,G_{p}(x_{p}) \bigr) =G ( \mathbf{x} ).
\end{eqnarray*}
This means that $\mathbf{X}\preceq_{\mathrm{lo}}\mathbf{Y}$. The other
part of
the theorem is proved similarly.
\end{pf*}
\begin{pf*}{Proof of Proposition~\ref{Prop-Smax}}
Let $\mathbf{X}$ and $\mathbf{Y}$ be independent MVNs with means
$\bolds{\mu}$ $\leq$ $\bolds{\nu}$ and common variance matrix
$\bolds{\Sigma}$. Clearly%
\[
\mathbb{P}\bigl(\mathbf{s}^{T}\mathbf{X}\leq\mathbf{s}^{T}%
\mathbf{Y}\bigr)=\Phi\biggl(-\frac{\mathbf{s}^{T}(\bolds{\mu}%
-\bolds{\nu})}{\sqrt{2\mathbf{s}^{T}\bolds{\Sigma} \mathbf{s}}}\biggr),
\]
where $\Phi$ is the DF of a standard normal RV. It follows that $\mathbb
{P}%
(\mathbf{s}^{T}\mathbf{X}\leq\mathbf{s}^{T}\mathbf{Y})$ is
maximized when the ratio $\mathbf{s}^{T}(\bolds{\nu-\mu}%
)/\sqrt{\mathbf{s}^{T}\bolds{\Sigma} \mathbf{s}}$ is maximized. From the
Cauchy--Schwarz inequality we have%
%
\begin{equation}\label{CS}
\frac{\mathbf{s}^{T}(\bolds{\nu}-\bolds{\mu})}{\sqrt{\mathbf{s}^{T}\bolds
{\Sigma} \mathbf{s}}}\leq\sqrt{(\bolds{\nu}-\bolds{\mu})^{T}\bolds{
\Sigma}^{-1}(\bolds{\nu}%
-\bolds{\mu})}
\end{equation}
for all $\mathbf{s}$. It is now easily verified that $\mathbf{s}
=\bolds{\Sigma}^{-1}(\bolds{\nu-\mu})$ maximizes the
left-hand side
of (\ref{CS}).
\end{pf*}
\begin{pf*}{Proof of Proposition~\ref{Prop-Max2}}
Let $Q_{q}$, $q=1,\ldots,4$, be the four quadrants. It is clear that maximizing
(\ref{Psi-nm}) is equivalent to maximizing%
%
\begin{equation}\label{Psi-nm-Pf}
\Psi_{n,m}^{\prime}(\mathbf{s})=\sum_{\mathbf{Z}_{ij}\in Q_{2}}
\mathbb{I}_{ ( \mathbf{s}^{T}\mathbf{Z}_{ij}\geq0 ) }
+\sum_{\mathbf{Z}_{ij}\in Q_{4}} \mathbb{I}_{ ( \mathbf{s}
^{T}\mathbf{Z}_{ij}\geq0 ) }.
\end{equation}
It is also clear that for any $\mathbf{s}$ the indicators $\mathbb
{I}%
_{ ( \mathbf{s}^{T}\mathbf{Z}_{ij}\geq0 ) }$ are
independent of the length of $\mathbf{Z}_{ij}$ which we therefore
take to
have length unity. Observe that the value of (\ref{Psi-nm-Pf}) is
constant in
the intervals $(\theta_{[ i]},\theta_{[ i+1]})$ where
$\theta_{[ i]}$ are defined in Algorithm~\ref{Alg-Max2}. At each point
$\theta_{[ i]}$, $i=0,\ldots,M+1$, the value of (\ref{Psi-nm-Pf}) may
increase or decrease. It follows that for all $\mathbf{s}\in
\mathcal{S}_{+}^{p-1}$ $\Psi_{n,m}^{\prime}(\mathbf{s})\in\{\Psi
_{n,m}^{\prime}(\mathbf{s}_{ [ 0 ] }),\ldots,\Psi_{n,m}%
^{\prime}(\mathbf{s}_{ [ M+1 ] })\}$ where $\mathbf{s}%
_{ [ i ] }$ are defined in Algorithm~\ref{Alg-Max2}. Therefore the
maximum value of (\ref{Psi-nm}) is an element of the above list. Now suppose
that $\mathbf{s}_{ [ i ] }$ is a global maximizer of
(\ref{Psi-nm-Pf}). Clearly either $\Psi_{n,m}^{\prime}(\mathbf
{s}_{ [
i ] })=\Psi_{n,m}^{\prime}(\mathbf{s}_{ [ i-1 ] })$ or
$\Psi_{n,m}^{\prime}(\mathbf{s}_{ [ i ] })=\Psi_{n,m}^{\prime
}(\mathbf{s}_{ [ i+1 ] })$ must hold, in which case any value
in $[\theta_{ [ i-1 ] },\theta_{ [ i ] }]$ or
$[\theta_{ [ i ] },\theta_{ [ i+1 ] }]$ is a
global maximizer. This concludes the proof.~%
\end{pf*}
\begin{pf*}{Proof of Theorem~\ref{Them-LST1}}
Using Hajek's projection and for any $\mathbf{s}$, we may write%
%
\begin{equation}\label{Hajek-P-1}\qquad\quad
\Psi_{n,m}(\mathbf{s})=\Psi(\mathbf{s})+n^{-1}\sum
_{i=1}^{n}\psi_{1} ( \mathbf{X}_{i},
\mathbf{s} ) +m^{-1}\sum_{j=1}%
^{m}
\psi_{2} ( \mathbf{Y}_{j},\mathbf{s} ) +R_{n,m} (
\mathbf{s} ),
\end{equation}
where%
\begin{eqnarray*}
\psi_{1} ( \mathbf{X}_{i},\mathbf{s} ) &=& \overline{G}
\bigl( \mathbf{s}^{T}\mathbf{X}_{i} \bigr) -\Psi(\mathbf
{s}%
),
\\
\psi_{2} ( \mathbf{Y}_{j},\mathbf{s} ) &=& F \bigl(
\mathbf{s}^{T}\mathbf{Y}_{j} \bigr) -\Psi(\mathbf{s})
\end{eqnarray*}
and $R_{n,m} ( \mathbf{s} ) $ is a remainder term. Here
$\overline{G} ( \mathbf{s}^{T}\mathbf{x} ) =\mathbb{P}%
(\mathbf{s}^{T}\mathbf{Y}\geq\mathbf{s}^{T}\mathbf{x})$,
$F ( \mathbf{s}^{T}\mathbf{y} ) =\mathbb{P}(\mathbf
{s}%
^{T}\mathbf{X}\leq\mathbf{s}^{T}\mathbf{y})$ and $\Psi
(\mathbf{s})=\mathbb{E}(\overline{G} ( \mathbf{s}^{T}%
\mathbf{X}_{i} ) )=\mathbb{E}(F ( \mathbf{s}^{T}%
\mathbf{Y}_{j} ) )$. Clearly $\mathbb{E}[\psi_{1} (
\mathbf{X}_{i},\mathbf{s} ) ]=\mathbb{E}[\psi_{2} (
\mathbf{Y}_{j},\mathbf{s} ) ]=0$ for all $i$ and $j$, so
by the
strong law of large numbers $n^{-1}\sum_{i=1}^{n}\psi_{1} (
\mathbf{X}_{i},\mathbf{s} ) $ and $m^{-1}\sum_{j=1}^{m}\psi_{2} ( \mathbf
{Y}_{j},\mathbf{s} ) $ both converge to zero
with probability one. Now,%
\begin{eqnarray*}
\sup_{\mathbf{s}\in\mathcal{S}_{+}^{p-1}}\bigl\llvert\Psi_{n,m}%
(\mathbf{s})-\Psi(
\mathbf{s})\bigr\rrvert&\leq&\sup_{\mathbf
{s}%
\in\mathcal{S}_{+}^{p-1}}\Biggl\llvert n^{-1}\sum
_{i=1}^{n}\psi_{1} (
\mathbf{X}_{i},\mathbf{s} ) \Biggr\rrvert+\sup_{\mathbf
{s}%
\in\mathcal{S}_{+}^{p-1}}\Biggl
\llvert m^{-1}\sum_{j=1}^{m}
\psi_{2} ( \mathbf{Y}_{j},\mathbf{s} ) \Biggr\rrvert\\
&&{}+
\sup_{\mathbf
{s}%
\in\mathcal{S}_{+}^{p-1}}\bigl\llvert R_{n,m} ( \mathbf{s} ) \bigr
\rrvert.
\end{eqnarray*}
The set $\mathcal{S}_{+}^{p-1}$ is compact, and the function $\psi_{1} (
\mathbf{x},\mathbf{s} ) $ is continuous in $\mathbf{s}%
\in\mathcal{S}_{+}^{p-1}$ for all values of $\mathbf{x}$ and
bounded [in
fact $\llvert\psi_{1} ( \mathbf{x},\mathbf{s} )
\rrvert\leq2$]. Thus the conditions in Theorem 3.1 in \citet{Das08}
are satisfied, and it follows that $\sup_{\mathbf{s}\in\mathcal{S}%
_{+}^{p-1}}\llvert n^{-1}\sum_{i=1}^{n}\psi_{1} ( \mathbf{X}%
_{i},\mathbf{s} ) \rrvert\stackrel{\mathrm{a.s.}}{\rightarrow}0$ as
$n\rightarrow\infty$. Similarly
$
\sup_{\mathbf{s}\in\mathcal
{S}_{+}^{p-1}%
}\llvert m^{-1}\*\sum_{i=1}^{m}\psi_{2} (
\mathbf{Y}_{i},\mathbf{s} )
\rrvert\stackrel{\mathrm{a.s.}}{\rightarrow}0
$
as $m\rightarrow\infty$.
Since\vspace*{1pt} $\Psi_{n,m}(\mathbf{s})$ is bounded all its moments exist;
therefore from Theorem 5.3.3 in \citet{Ser80} we have that with
probability one $R_{n,m} ( \mathbf{s} ) =o (
1/N )$. Moreover it is clear that the latter holds uniformly for
all~$\mathbf{s}$. Thus,
\[
\sup_{\mathbf{s}\in\mathcal{S}_{+}^{p-1}}\bigl\llvert\Psi_{n,m}%
(\mathbf{s})-\Psi(
\mathbf{s})\bigr\rrvert\stackrel{\mathrm{a.s.}} {\rightarrow}0\qquad
\mbox{as }n,m\rightarrow
\infty.
\]
By assumption $\Psi(\mathbf{s}_{\max})>\Psi(\mathbf{s})$ for all
$\mathbf{s}\in\mathcal{S}_{+}^{p-1}\setminus\mathbf{s}_{\max}$
so we
can apply Theorem~2.12 in \citet{Kos08} to conclude that%
\[
\widehat{\mathbf{s}}_{\max}\stackrel{\mathrm{a.s.}} {\rightarrow}\mathbf
{s}%
_{\max};%
\]
that is, $\widehat{\mathbf{s}}_{\max}$ is strongly consistent. This
completes
the first part of the proof.

Since the densities of $\mathbf{X}$ and $\mathbf{Y}$ are
differentiable, it follows that $\Psi(\mathbf{s})$ is continuous
and twice
differentiable. In particular at $\mathbf{s}_{\max}\in\mathcal{S}%
_{+}^{p-1}$, the matrix $-\nabla^{2}\Psi(\mathbf{s}_{\max})$ exists
and is
positive definite. A Taylor expansion implies that%
\[
\sup_{\llVert\mathbf{s}-\mathbf{s}_{\max}\rrVert<\delta
}%
\Psi(\mathbf{s})-\Psi(\mathbf{s}_{\max})
\leq-C\delta^{2}.
\]
It is also obvious that%
\[
n^{-1}\sum_{i=1}^{n}
\psi_{1} ( \mathbf{X}_{i},\mathbf{s}_{\max
} )
=O_{p}(1/\sqrt{N})
\]
and
\[
m^{-1}\sum
_{j=1}^{m}\psi_{2} ( \mathbf{Y}_{j},
\mathbf{s}_{\max} ) =O_{p}(1/\sqrt{N}).
\]

Finally as noted above $\llvert R_{n,m} ( \mathbf{s} )
\rrvert=O(1/N)$ for all $\mathbf{s}$ as $n,m\rightarrow\infty$.
Therefore by Theorem 1 in \citet{She93} we have that%
%
\begin{equation}\label{root-n-rate}
\widehat{\mathbf{s}}_{\max}=\mathbf{s}_{\max}+O_{p}
\bigl(N^{-1/2}%
\bigr);%
\end{equation}
that is, $\widehat{\mathbf{s}}_{\max}$ converges to $\mathbf
{s}_{\max} $
at a $N^{1/2}$ rate. This completes the second part of the proof.

The functions $\Psi(\mathbf{s}),\psi_{1} ( \mathbf{X}%
_{i},\mathbf{s} ) $ and $\psi_{2} ( \mathbf{Y}%
_{j},\mathbf{s} ) $ on the right-hand side of (\ref
{Hajek-P-1}) all
admit a quadratic expansion. A bit of algebra shows that for
$\mathbf{s}$
in an $O_{p}(N^{-1/2})$ neighborhood of $\mathbf{s}_{\max}$, we have
%
\begin{eqnarray}\label{Hajek-P-2}
\Psi_{n,m}(\mathbf{s})&=&\Psi(\mathbf{s}_{\max})+ (
\mathbf{s}%
-\mathbf{s}_{\max} )^{T}
\frac{\mathbf{M}_{n,m}}{N^{1/2}}\nonumber\\[-8pt]\\[-8pt]
&&{}+\frac{1}{2} ( \mathbf{s}-\mathbf{s}_{\max}
)^{T}\mathbf{V} ( \mathbf{s}-\mathbf{s}_{\max} ) +o_{p}
( 1/N ),\nonumber
\end{eqnarray}
where%
\[
\mathbf{M}_{n,m}=\frac{1}{\sqrt{\lambda_{n,m}}}\frac{\sum_{i=1}^{n}%
\nabla\psi_{1} ( \mathbf{X}_{i},\mathbf{s}_{\max} )
}{n^{1/2}}+\frac{1}{\sqrt{1-\lambda_{n,m}}}
\frac{\sum_{j=1}^{m}\nabla
\psi_{2} ( \mathbf{Y}_{j},\mathbf{s}_{\max} ) }{m^{1/2}}%
\]
$\lambda_{n,m}=n/N$, for $j=1,2$ the function $\nabla\psi_{j} (
\cdot,\mathbf{s}_{\max} )^{T}$ is the gradient of $\psi_{j} (
\cdot,\mathbf{s} ) $ evaluated at $\mathbf{s}_{\max}$, and the
matrix $\mathbf{V}$ is given by%
\[
\mathbf{V}=\mathbb{E} \bigl( \nabla^{2}\psi_{1} ( \mathbf{X},
\mathbf{s}_{\max} ) \bigr) +\mathbb{E} \bigl( \nabla^{2}
\psi_{2} ( \mathbf{Y},\mathbf{s}_{\max} ) \bigr).
\]
Note that the $o_{p} ( 1/N ) $ term in (\ref{Hajek-P-2}) absorbs
$R_{n,m} ( \mathbf{s} ) $ in (\ref{Hajek-P-1}) as well as the
higher-order terms in the quadratic expansions of $\Psi(\mathbf{s}%
),n^{-1}\sum_{i=1}^{n}\psi_{1} ( \mathbf{X}_{i},\mathbf{s}%
) $ and $m^{-1}\sum_{j=1}^{m}\psi_{2} ( \mathbf{Y}%
_{j},\mathbf{s} ) $. Now by the CLT and Slutzky's theorem, we have
that%
\[
\mathbf{M}_{n,m}\Rightarrow N ( 0,\bolds{\Delta} ),
\]
where
\begin{eqnarray*}
\bolds{\Delta}&=&\frac{1}{\lambda}\mathbb{E} \bigl( \nabla\psi_{1} (
\mathbf{X},\mathbf{s}_{\max} ) \nabla\psi_{1} ( \mathbf{X},
\mathbf{s}_{\max} )^{T} \bigr) \\
&&{}+\frac
{1}{1-\lambda
}\mathbb{E}
\bigl( \nabla\psi_{2} ( \mathbf{Y},\mathbf{s}_{\max
} ) \nabla
\psi_{2} ( \mathbf{Y},\mathbf{s}_{\max} )^{T} \bigr).
\end{eqnarray*}
Finally it follows by Theorem 2 in \citet{She93} that%
\[
N^{1/2} ( \widehat{\mathbf{s}}_{\max}-\mathbf{s}_{\max
} )
\Rightarrow N ( 0,\bolds{\Sigma} ),
\]
where $\bolds{\Sigma}=\mathbf{V}^{-1}\bolds{\Delta V}^{-1}$,
completing the proof.
\end{pf*}
\begin{pf*}{Proof of Theorem~\ref{Them-LST2}}
Suppose that $\mathbf{X}$ and $\mathbf{Y}$ are
discrete RVs with finite support. Let $p_{a}=\mathbb{P}(\mathbf
{X}%
=\mathbf{x}_{a})>0$ and $q_{b}=\mathbb{P}(\mathbf
{Y}=\mathbf{y}%
_{b})>0$ where $a=1,\ldots,A$ and $b=1,\ldots,B;$ $A$ and $B$ are finite.
Define the set $\mathbf{S}_{ab}=\{\mathbf{s}\in\mathcal{S}_{+}%
^{p-1}\dvtx \mathbf{s}^{T}\mathbf{x}_{a}\leq\mathbf{s}^{T}%
\mathbf{y}_{b}\}$. A simple argument shows that%
\[
\Psi(\mathbf{s})=\mathbb{P}\bigl(\mathbf{s}^{T}\mathbf{X}%
\leq
\mathbf{s}^{T}\mathbf{Y}\bigr)=\sum_{a=1}^{A}
\sum_{b=1}^{B}%
\mathbb{I}_{(\mathbf{s}\in\mathbf{S}_{ab})}%
p_{a}q_{b}=\sum
_{k=1}^{K}\alpha_{k}
\mathbb{I}_{(\mathbf{s}%
\in\mathbf{S}_{k})},
\]
where $K\leq2^{AB}$ is finite, the sets $\mathbf{S}_{k}$ are distinct,
$\bigcup_{k=1}^{K}\mathbf{S}_{k}=\mathcal{S}_{+}^{p-1}$ and $\alpha_{k}
=\sum_{(a,b)\in J_{k}}p_{a}q_{b}$ with $J_{k}= \{ (
a,b )\dvtx \mathbf{S}_{k}\cap\mathbf{S}_{ab}\neq\varnothing
\}
$. Thus $\Psi(\mathbf{s})$ is a simple function on $\mathcal{S}_{+}^{p-1}$,
and $\mathbf{S}_{\max}$ is the set associated with the largest
$\alpha_{k}$. We will assume, without any loss of generality, that
$\alpha_{1}%
>\alpha_{k}$ for all $k\geq2$. Now note that%
\[
\Psi_{n,m}(\mathbf{s})=\frac{1}{nm}\sum
_{i=1}^{n}\sum_{j=1}^{m}%
\mathbb{I}_{ ( \mathbf{s}^{T}\mathbf{X}_{i}\leq\mathbf
{s}%
^{T}\mathbf{Y}_{j} ) }=\sum_{a=1}^{A}
\sum_{b=1}^{B}\frac
{n_{a}m_{b}%
}{nm}
\mathbb{I}_{(\mathbf{s}\in\mathbf{S}_{ab}
)},
\]
where\vspace*{1pt} $n_{a}=\sum_{i=1}^{n}\mathbb{I}_{ ( \mathbf{X}%
_{i}=\mathbf{x}_{a} ) },m_{b}=\sum_{j=1}^{m}\mathbb
{I}_{ (
\mathbf{Y}_{i}=\mathbf{y}_{b} ) }$ where $a=1,\ldots,A$ and
$b=1,\ldots,B$. Clearly $\Psi_{n,m}(\mathbf{s})$ is also a simple
function. Moreover for large enough $n$ and $m$ we will have $n_{a}>0$ and
$m_{b}>0$ for all $a=1,\ldots,A$ and $b=1,\ldots,B$, and consequently
$\Psi_{n,m}(\mathbf{s})$ is defined over the same sets as $\Psi
(\mathbf{s})$, that is,
\[
\Psi_{n,m}(\mathbf{s})=\sum_{k=1}^{K}
\widehat{\alpha}_{k}\mathbb{I}%
_{(\mathbf{s}\in\mathbf{S}_{k})},
\]
where $\widehat{\alpha}_{k}=\sum_{(a,b)\in J_{k}}\widehat{p}_{a}\widehat
{q}_{b}$ with $\widehat{p}_{a}=n_{a}/n$ and $\widehat{q}_{b}=m_{b}/m$.
Furthermore the maximizer of $\Psi_{n,m}(\mathbf{s})$ is any
$\mathbf{s}\in\mathbf{S}_{k}$ provided that $\mathbf
{S}_{k}$ is
associated with the largest~$\widehat{\alpha}_{k}$. Hence,%
%
\begin{eqnarray}\label{bound-1}
\mathbb{P}(\widehat{\mathbf{s}}_{\max}
\notin\mathbf{S}_{\max})
&=&
\mathbb{P}\Bigl(\arg\max_{1\leq k\leq K}\widehat{\alpha}_{k}\neq1
\Bigr)=\mathbb{P}%
\Biggl(\bigcup_{k=2}^{K}
\{ \widehat{\alpha}_{1}\leq\widehat{\alpha}_{k} \} \Biggr)
\nonumber\\[-8pt]\\[-8pt]
&\leq& \sum_{k=2}^{K}\mathbb{P}(\widehat{
\alpha}_{1}\leq\widehat{\alpha}%
_{k})\leq(K-1)
\max_{2\leq k\leq K}\mathbb{P}(\widehat{\alpha}_{1}\leq\widehat{
\alpha}_{k}).\nonumber
\end{eqnarray}
A bit of rearranging shows that%
\begin{eqnarray*}
\mathbb{P}(\widehat{\alpha}_{1}\leq\widehat{\alpha}_{k})&=&
\mathbb{P}%
\biggl(\sum_{(a,b)\in J_{1}}
\widehat{p}_{a}\widehat{q}_{b}-\sum
_{(a,b)\in
J_{k}%
}\widehat{p}_{a}\widehat{q}_{b}
\leq0\biggr)\\
&=&\mathbb{P}\Biggl(n^{-1}m^{-1}\sum
_{i=1}%
^{n}\sum
_{j=1}^{m}Z_{ij}^{ ( k ) }\leq0\Biggr),
\end{eqnarray*}
where
\[
Z_{ij}^{ ( k ) }=\sum_{(a,b)\in J_{1}\setminus J_{k}}%
\mathbb{I}_{ ( \mathbf{X}_{i}=\mathbf{x}_{a} ) }%
\mathbb{I}_{ ( \mathbf{Y}_{j}=\mathbf{y}_{b} ) }%
-\sum
_{(a,b)\in J_{k}\setminus J_{1}}\mathbb{I}_{ ( \mathbf{X}%
_{i}=\mathbf{x}_{a} ) }\mathbb{I}_{ ( \mathbf{Y}%
_{j}=\mathbf{y}_{b} ) }.
\]
Note that $Z_{ij}^{ ( k ) }$ may be viewed as a kernel of a two
sample $U$-statistic. Moreover%
\[
-\llvert J_{k}\setminus J_{1}\rrvert\leq
Z_{ij}^{ ( k )
}\leq\llvert J_{k}\setminus
J_{1}\rrvert
\]
is bounded (here $\llvert\cdot\rrvert$ denotes set
cardinality) and
$\mathbb{E}(Z_{ij}^{ ( k ) })=\mu_{k}=\mathbb{E} (
\widehat{\alpha}_{1}-\widehat{\alpha}_{k} ) >0$ by assumption. Applying\vadjust{\goodbreak}
Theorem 2 and the derivations in Section~5b in \citet{Hoe63} we have
that%
%
\begin{equation}\label{bound-2}
\mathbb{P}\Biggl(n^{-1}m^{-1}\sum
_{i=1}^{n}\sum_{j=1}^{m}Z_{ij}^{ (
k )
}
\leq0\Biggr)\leq\exp\biggl( -\frac{N\min( \lambda_{n,m},1-\lambda
_{n,m} ) \mu_{k}^{2}}{2 ( \llvert J_{k}\setminus J_{1}%
\rrvert+\llvert J_{1}\setminus J_{k}\rrvert)^{2}%
} \biggr),\hspace*{-35pt}
\end{equation}
where $\lambda_{n,m}=n/N\rightarrow\lambda\in( 0,1 ) $ as
$n,m\rightarrow\infty$. Finally from (\ref{bound-1}) and (\ref
{bound-2}) we
have that%
\[
\mathbb{P}(\widehat{\mathbf{s}}_{\max}\notin\mathbf{S}_{\max
})\leq
C_{1}\exp( -C_{2}N ),
\]
where $C_{1}=K-1$ and $C_{2}=\min( \lambda,1-\lambda) \min$
$2\mu_{k}^{2} ( \llvert J_{k}\setminus J_{1}\rrvert+\llvert
J_{1}\setminus J_{K}\rrvert)^{-2}$ completing the proof.
\end{pf*}
\begin{pf*}{Proof of Theorem~\ref{Them-LST3}}
Choose $\varepsilon>0$. We have already seen that under
the stated conditions, $\Psi(\mathbf{s})$ is continuous, and
therefore for
each $\mathbf{s}$ the set $B_{\mathbf{s},\varepsilon}= \{
\mathbf{s}^{\prime}\dvtx \llvert\Psi(\mathbf{s}^{\prime})-\Psi
(\mathbf{s})\rrvert<\varepsilon\} $ is open. The collection
$\{B_{\mathbf{s},\varepsilon}\dvtx \mathbf{s}\in\mathcal
{S}_{+}^{p-1}\}$ is
an open cover for $\mathcal{S}_{+}^{p-1}$. Since $\mathcal
{S}_{+}^{p-1}$ is
compact there exists a finite subcover $B_{\mathbf
{s}_{1},\varepsilon
},\ldots,B_{\mathbf{s}_{_{K}},\varepsilon}$ for $\mathcal{S}_{+}^{p-1}$
where $K<\infty$. Hence each $\mathbf{s}$ belongs to some
$B_{\mathbf{s}_{_{i}},\varepsilon}$ and therefore%
\begin{eqnarray*}
&&
\sup_{\mathbf{s}\in\mathcal{S}_{+}^{p-1}}\bigl\llvert\Psi_{N}%
(\mathbf{s})-\Psi(
\mathbf{s})\bigr\rrvert\\
&&\qquad\leq
\max_{1\leq i\leq
K}\sup_{\mathbf{s}\in B_{\mathbf{s}_{i},\varepsilon}}\bigl\llvert
\Psi_{N}(\mathbf{s})-\Psi_{N}(\mathbf{s}_{i})\bigr
\rrvert
\\
&&\qquad\quad{} +\max_{1\leq i\leq K}\bigl\llvert\Psi_{N}(\mathbf{s}_{i})-
\Psi(\mathbf{s}_{i})\bigr\rrvert+\max_{1\leq i\leq K}
\sup_{\mathbf
{s}\in
B_{\mathbf{s}_{i},\varepsilon}}\bigl\llvert\Psi(\mathbf{s}_{i}%
)-\Psi(
\mathbf{s})\bigr\rrvert.
\end{eqnarray*}
By construction $\sup_{\mathbf{s}\in B_{\mathbf
{s}_{i},\varepsilon}%
}\llvert\Psi(\mathbf{s}_{i})-\Psi(\mathbf{s})\rrvert
<\varepsilon$ for all $i$. By the law of large numbers $\Psi_{N}%
(\mathbf{s}_{i})\stackrel{\mathrm{a.s.}}{\rightarrow}\Psi(\mathbf
{s}_{i})$ as
$N\rightarrow\infty$ for each $i=1,\ldots,K$. Since $K$ is finite,
\[
\max_{1\leq i\leq K}\bigl\llvert\Psi_{N}(\mathbf{s}_{i})-
\Psi(\mathbf{s}_{i})\bigr\rrvert\leq\sum_{i=1}^{K}
\bigl\llvert\Psi_{N}(\mathbf{s}_{i})-\Psi(
\mathbf{s}_{i})\bigr\rrvert\stackrel{\mathrm{a.s.}} {\rightarrow}0.
\]
Now for any $\mathbf{s}^{\prime}\in B_{\mathbf
{s}_{_{i}},\varepsilon
}$, $\Psi_{N}(\mathbf{s}^{\prime})\stackrel{\mathrm{a.s.}}{\rightarrow}%
\Psi(\mathbf{s}^{\prime})$, $\Psi_{N}(\mathbf{s}_{i})\stackrel
{\mathrm{a.s.}}{\rightarrow}\Psi(\mathbf{s}_{i})$ and $\llvert\Psi
(\mathbf{s}^{\prime})-\Psi(\mathbf{s}_{i})\rrvert
<\varepsilon$.
This implies that we can choose $N$ large enough so $\llvert\Psi
_{N}(\mathbf{s}^{\prime})-\Psi_{N}(\mathbf{s}_{i})\rrvert
<2\varepsilon$. Moreover this bound holds for all $\mathbf
{s}^{\prime}\in
B_{\mathbf{s}_{_{i}},\varepsilon}$ and $i$ so
\[
\lim\sup\bigl\llvert\Psi_{N}(\mathbf{s})-\Psi(\mathbf{s})\bigr\rrvert
\leq3\varepsilon
\]
on $\mathcal{S}_{+}^{p-1}$. Since $\varepsilon$ is arbitrary we
conclude that
$\sup_{\mathbf{s}\in\mathcal{S}_{+}^{p-1}}\llvert\Psi_{N}%
(\mathbf{s})-\Psi(\mathbf{s})\rrvert\stackrel
{\mathrm{a.s.}}{\rightarrow}0$
as $N\rightarrow\infty$. By assumption $\Psi(\mathbf{s}_{\max}%
)>\Psi(\mathbf{s})$ for all $\mathbf{s}\in\mathcal{S}_{+}%
^{p-1}\setminus\mathbf{s}_{\max}$, so we can apply Theorem 2.12 in \citet{Kos08} to conclude that%
\[
\widehat{\mathbf{s}}_{\max}\stackrel{\mathrm{a.s.}} {\rightarrow}\mathbf
{s}%
_{\max},%
\]
that is, $\widehat{\mathbf{s}}_{\max}$ is strongly consistent. This
completes
the first part of the proof.\vadjust{\goodbreak}

We have already seen that%
\[
\sup_{\llVert\mathbf{s}-\mathbf{s}_{\max}\rrVert<\delta
}%
\Psi(\mathbf{s})-\Psi(\mathbf{s}_{\max})
\leq-C\delta^{2}%
\]
holds. We now need to bound $\mathbb{E}^{\ast}\sup_{\llVert
\mathbf{s}%
-\mathbf{s}_{\max}\rrVert<\delta}N^{1/2}\llvert(\mathbb{P}%
_{N}(\overline{\Psi} ( \mathbf{s} ) -\overline{\Psi} (
\mathbf{s}_{\max} ) )\rrvert$, where $\mathbb{E}^{\ast}$
denotes the outer expectation and $\overline{\Psi} ( \mathbf{s}%
) = \mathbb{I}_{ ( \mathbf{s}^{T}\mathbf{Z}\geq
0 )
}-\Psi(\mathbf{s})$. We first note that the bracketing entropy of the
upper half-planes is of the order $\delta/\varepsilon^{2}$. The envelope
function of the class $\mathbb{I}_{ ( \mathbf{s}^{T}\mathbf
{z}%
\geq0 ) }-\mathbb{I}_{ ( \mathbf{s}_{\max}^{T}\mathbf
{z}%
\geq0 ) }$ where $\llVert\mathbf{s}-\mathbf{s}_{\max
}\rrVert<\delta$ is bounded by $\mathbb{I}_{ ( \mathbf{s}%
^{T}\mathbf{z}\geq0>\mathbf{s}_{\max}^{T}\mathbf{z} )
}+\mathbb{I}_{ ( \mathbf{s}_{\max}^{T}\mathbf{z}\geq
0>\mathbf{s}^{T}\mathbf{z} ) }$ whose squared $L_{2}$ norm
is%
%
\begin{equation}
\label{vv-2-2} \mathbb{P} \bigl( \mathbf{s}^{T}\mathbf{Z}\geq0>\mathbf
{s}_{\max
}^{T}\mathbf{Z} \bigr) +\mathbb{P} \bigl(
\mathbf{s}_{\max}%
^{T}\mathbf{Z}\geq0>
\mathbf{s}^{T}\mathbf{Z} \bigr).
\end{equation}
Note that we may replace the RV $\mathbf{Z}$ in (\ref{vv-2-2}) with
the RV
$\mathbf{Z}^{\prime}=\mathbf{Z}/\llVert\mathbf{Z}\rrVert$
whose mass is concentrated on the unit sphere. The condition that
$\llVert
\mathbf{s}-\mathbf{s}_{\max}\rrVert<\delta$ implies that the
angle between $\mathbf{s}$ and $\mathbf{s}_{\max}$ is of the order
$O ( \delta), $ and therefore $\mathbb{P} ( \mathbf{s}
^{T}\mathbf{Z}^{\prime}\geq0>\mathbf{s}_{\max}^{T}\mathbf{Z}
^{\prime} ) $ is computed as surface integral on a spherical wedge with
maximum width~$\delta$. It follows that (\ref{vv-2-2}) is bounded by
$2A_{p-1}\delta\llVert h^{\prime}\rrVert_{\infty}$ where
$A_{p-1}$ is
the area of $\mathcal{S}_{+}^{p-1}$, and $\llVert h^{\prime}\rrVert
_{\infty}$ is the supremum of the density of $\mathbf{Z}^{\prime}$.
Clearly $\llVert h^{\prime}\rrVert_{\infty}<\infty$ since the density
of $\mathbf{Z}$ is bounded by assumption. Thus by Corollary 19.35
in \citet{van00} we have%
\[
\mathbb{E}^{\ast}\sup_{\llVert\mathbf{s}-\mathbf{s}_{\max
}\rrVert<\delta}N^{1/2}\bigl\llvert
\mathbb{P}_{N}\bigl(\overline{\Psi} ( \mathbf{s} ) -\overline{\Psi} (
\mathbf{s}_{\max} ) \bigr)\bigr\rrvert\leq C\delta^{1/2}.
\]
It now follows that%
%
\begin{equation}\label{KP-1}
\Psi_{N}(\widehat{\mathbf{s}}_{\max})\geq\sup_{\mathbf{s}%
\in\mathcal{S}_{+}^{p-1}}
\Psi_{N}(\mathbf{s})-o_{p}\bigl(N^{-2/3}\bigr),
\end{equation}
which implies by Theorem 5.52 in \citet{van00} and Theorem 14.4 of
\citet{Kos08} that%
\[
\widehat{\mathbf{s}}_{\max}=\mathbf{s}_{\max}+O_{p}
\bigl( N^{-1/3} \bigr);
\]
that is, $\widehat{\mathbf{s}}_{\max}$ converges to $\mathbf
{s}_{\max} $
at a cube root rate. This completes the second part of the proof.

The limit distribution is derived by verifying the conditions in
Theorem 1.1 of \citet{KimPol90}, denoted henceforth by KP. First
note that (\ref{KP-1}) is condition (i) in KP. Since
$\widehat{\mathbf{s}}_{\max}$ is consistent, condition (ii) also holds,
and condition (iii) holds by assumption. The differentiability of the
density of $\mathbf{Z}$ implies that $\Psi(\mathbf{s})$ is twice
differentiable. The uniqueness of the maximizer implies that
$-\nabla^{2}\Psi(\mathbf{s}_{\max})$ is positive definite, and hence
condition (iv) holds; see also Example 6.4 in KP for related
calculations. Condition (v) in KP is equivalent to the existence of the
limit $H ( \mathbf{u},\mathbf{v} )
=\lim_{\alpha\rightarrow\infty}\alpha\mathbb {E} ( \overline{\Psi} (
\mathbf{s}_{\max}+\mathbf{u}/\alpha) \overline{\Psi} (
\mathbf{s}_{\max}+\mathbf{v}/\alpha)
) $ which can be\vadjust{\goodbreak} rewritten as%
\begin{eqnarray*}
&&\lim_{\beta\rightarrow0}\frac{1}{\beta} \bigl[ \mathbb{P} \bigl(( \mathbf
{s}_{\max}+
\beta\mathbf{u})^{T}\mathbf{Z}\geq0,(\mathbf{s}_{\max}+\beta
\mathbf{v})^{T}\mathbf{Z}\geq0 \bigr) \\
&&\qquad\hspace*{11pt}{}-\mathbb{P} \bigl( (
\mathbf{s}_{\max}+\beta\mathbf{u})^{T}%
\mathbf{Z}\geq0
\bigr) \mathbb{P} \bigl( (\mathbf{s}_{\max}%
+\beta
\mathbf{v})^{T}\mathbf{Z}\geq0 \bigr) \bigr].
\end{eqnarray*}
With some algebra we find that this limit exists and equals%
\begin{eqnarray*}
&&
\int_{\mathbf{z}\in\mathcal{S}_{+}^{p-1}}\delta\bigl( \mathbf{s} 
_{\max}^{T}
\mathbf{z} \bigr) \bigl(\mathbf{z}^{T}(\mathbf{u}+\mathbf{v}%
)\bigr)h (
\mathbf{z} ) \,d\mathbf{z}\\
&&\qquad{}-\int_{\mathbf{z}%
\in\mathcal{S}_{+}^{p-1}}\delta\bigl(
\mathbf{s}_{\max}^{T}%
\mathbf{z} \bigr) \bigl(
\mathbf{z}^{T}\mathbf{u}\bigr)h ( \mathbf{z} ) \,d\mathbf{z}\int
_{\mathbf{z}\in\mathcal{S}%
_{+}^{p-1}}\delta\bigl( \mathbf{s}_{\max}^{T}
\mathbf{z} \bigr) \bigl(\mathbf{z}^{T}\mathbf{v}\bigr)h ( \mathbf{z} ) \,d
\mathbf{z},
\end{eqnarray*}
where $\delta( \mathbf{s}_{\max}^{T}\mathbf{z} ) $ is the usual Dirac
function; hence integration is with respect to the surface measure on $
\{ \mathbf{s}_{\max}^{T}\mathbf{z}=0 \} $. It follows that condition
also (v) holds. Conditions (vi) and (vii) were verified in the second part
of the proof. Thus we
may apply Theorem 1.1 in KP to get%
\[
N^{1/3}(\widehat{\mathbf{s}}_{\max}-\mathbf{s}_{\max
})
\Rightarrow\arg\max\bigl\{ -Q ( \mathbf{s} ) +\mathbb{W}(\mathbf{s}%
)\dvtx
\mathbf{s}\in\mathcal{S}_{+}^{p-1} \bigr\},
\]
where by KP $Q ( \mathbf{s} ) =\mathbf{s}^{T}\nabla^{2}
\Psi(\mathbf{s}_{\max})\mathbf{s}$ and $\mathbb{W}(\mathbf{s})$ is a zero
mean Gaussian process with covariance function $H (
\mathbf{u},\mathbf{v} )$. This completes the proof.
\end{pf*}
\begin{pf*}{Proof of Proposition~\ref{Prop-Unique}}
Note that%
\begin{eqnarray*}
\Psi(\mathbf{s})&=&\mathbb{P}\bigl(\mathbf{s}^{T}\mathbf{X}%
\leq
\mathbf{s}^{T}\mathbf{Y}\bigr)=\mathbb{P} \bigl( \mathbf{s}%
^{T}
\mathbf{Z}\geq0 \bigr) =\mathbb{P}\biggl(\frac{\mathbf{s}%
^{T}\mathbf{Z}-\mathbf{s}^{T}\bolds{\delta}}{\sqrt{\mathbf{s}%
^{T}\bolds{\Sigma} \mathbf{s}}}\geq-\frac{\mathbf{s}^{T}\bolds
{\delta}%
}{\sqrt{\mathbf{s}^{T}\bolds{\Sigma} \mathbf{s}}}
\biggr)\\
&=&1-F\biggl(-\frac
{\mathbf{s}%
^{T}\bolds{\delta}}{\sqrt{\mathbf{s}^{T}\bolds{\Sigma} \mathbf{s}}}\biggr).
\end{eqnarray*}
Now, by assumption the DF $F$ is independent of $\mathbf{s}$. Therefore
$\Psi(\mathbf{s})$ is uniquely maximized on $\mathcal{S}_{+}^{p-1}$
if and
only if the function%
\[
\varkappa( \mathbf{s} ) =\frac{\mathbf{s}^{T}%
\bolds{\delta}}{\sqrt{\mathbf{s}^{T}\bolds{\Sigma} \mathbf{s}}}%
\]
is uniquely maximized on $\mathcal{S}_{+}^{p-1}$. If $\bolds{\Sigma=I}$,
then $\varkappa( \mathbf{s} ) =\mathbf{s}^{T}%
\bolds{\delta}$, and we wish to maximize a linear function on
$\mathcal{S}_{+}^{p-1}$. It is easily verified (by using ideas from linear
programming) that the maximizer is unique if $\bolds{\delta}%
\geq\mathbf{0}$ which is true by assumption. Incidentally, it is
easy to
show directly that $\varkappa( \mathbf{s} ) $ is
maximized at
$\mathbf{s}^{\ast}/\llVert\mathbf{s}^{\ast}\rrVert$
where%
\[
\mathbf{s}^{\ast}= ( \delta_{1}\mathbb{I}_{(\delta_{1}%
\geq0)},\ldots,\delta_{p}\mathbb{I}_{(\delta_{p}\geq0)} ).
\]
Now let $\bolds{\Sigma}\neq\mathbf{I}$ and assume that a unique
maximizer does not exist; that is, suppose that $\varkappa(
\mathbf{s}%
) $ is maximized by both $\mathbf{s}_{1}$ and $\mathbf{s}_{2}$.
It is clear that $\varkappa( \lambda_{1}\mathbf{s}_{1} )
=\varkappa( \lambda_{2}\mathbf{s}_{2} ) $ for all $\lambda_{1},\lambda
_{2}>0;$ that is, the value of $\varkappa( \cdot
) $ is
constant along rays through the origin. The rays passing\vadjust{\goodbreak} through
$\mathbf{s}_{1}$ and $\mathbf{s}_{2}$, respectively, intersect the
ellipsoid $\mathbf{s}^{T}\bolds{\Sigma} \mathbf{s}=1$ at the points
$\mathbf{p}_{1}$ and $\mathbf{p}_{2}$. It follows that
$\varkappa
( \mathbf{p}_{1} ) =\varkappa( \mathbf{p}%
_{2} )$, moreover $\mathbf{p}_{1}$ and $\mathbf{p}_{2}$ maximize
$\varkappa( \cdot) $ on the ellipsoid. Now since
$\mathbf{p}_{1}^{T}\bolds{\Sigma} \mathbf{p}_{1}=1=\mathbf{p}_{2}%
^{T}\bolds{\Sigma} \mathbf{p}_{2}$ we must have $\mathbf{p}_{1}^{T}%
\bolds{\delta}=\mathbf{p}_{2}^{T}\bolds{\delta}$.
Recall that
a linear function on ellipsoid is uniquely maximized (just like on a sphere;
see the comment above). Therefore we must have $\mathbf{p}%
_{1}=\mathbf{p}_{2}$ which implies that $\mathbf{s}_{1}=
\mathbf{s}_{2}$ as required.
\end{pf*}
\begin{pf*}{Proof of Theorem~\ref{Thm-TestStat}}
If $\mathbf{X}=_{\st }\mathbf{Y}$, then for all $\mathbf{s}$
we have
$\mathbf{s}^{T}\mathbf{X}=_{\st }\mathbf{s}^{T}\mathbf
{Y}$. By
assumption both $\mathbf{s}^{T}\mathbf{X}$ and $\mathbf{s}%
^{T}\mathbf{Y}$ are continuous RVs, so $\mathbb{P} ( \mathbf
{s}%
^{T}\mathbf{X}\leq\mathbf{s}^{T}\mathbf{Y} ) =1/2$. Suppose
now that both $\mathbf{X}\preceq_{\st }\mathbf{Y}$ and $\mathbb{P}
( \mathbf{s}^{T}\mathbf{X}\leq\mathbf{s}^{T}\mathbf
{Y}%
) >1/2$ for some $\mathbf{s}\in\mathcal{S}_{+}^{p-1} $, hold. Then
we must have $\mathbf{X}\prec_{\lst }\mathbf{Y}$. Since
$\mathbf{X}%
\preceq_{\st }\mathbf{Y}$ we have $X_{j}\preceq_{\st }Y_{j}$ for $1\leq
j\leq
p$. One of these inequalities must be strict; otherwise $\mathbf{X}%
=_{\st }\mathbf{Y}$ contradicts the fact that $\mathbf{X}\prec_{\lst }\mathbf
{Y}$. Now use Theorem 1 in \citet{DavPed11} to
complete the proof.
\end{pf*}
\begin{pf*}{Proof of Theorem~\ref{Thm-TestLimit}}
The functions $\psi_{1}$ and $\psi_{2}$ defined in the proof of Theorem
\ref{Them-LST1} are Donsker; cf. Example 19.7 in \citet{van00}. Hence
by the theory of empirical processes applied to (\ref{Hajek-P-1}), we
find that%
%
\begin{equation}\label{KP-2}
N^{1/2}\bigl(\Psi_{n,m}(\mathbf{s})-\Psi(\mathbf{s})\bigr)
\Rightarrow\mathbb{G}(\mathbf{s}),
\end{equation}
where $\mathbb{G}(\mathbf{s})$ is a zero mean Gaussian process, and
convergence holds for all $\mathbf{s}\in\mathcal{S}_{+}^{p-1}$. We also
note that (\ref{KP-2}) is a two-sample $U$-processes. A central limit theorem
for such processes is described by \citet{Neu04}. Hence by the continuos
mapping theorem, and under $H_{0}$, we have $N^{1/2}(\Psi_{n,m}(\widehat
{\mathbf{s}}_{\max})-1/2))\Rightarrow\sup_{\mathbf{s}\in
\mathcal{S}_{+}^{p-1}}\mathbb{G}(\mathbf{s})$ where the covariance
function of $\mathbb{G}(\mathbf{s})$, denoted by $C (
\mathbf{u},\mathbf{v} ) $, is given by%
%
\begin{eqnarray}\label{Cu,v-H0}
&&
\frac{1}{\lambda}\mathbb{P}\bigl(\mathbf{u}^{T}\mathbf{X}_{1}%
\leq\mathbf{u}^{T}\mathbf{X}_{2},\mathbf{v}^{T}\mathbf
{X}%
_{1}\leq\mathbf{v}^{T}\mathbf{X}_{3}
\bigr)\nonumber\\[-8pt]\\[-8pt]
&&\qquad{}+\frac{1}{1-\lambda}%
\mathbb{P}\bigl(\mathbf{u}^{T}
\mathbf{X}_{1}\leq\mathbf{u}%
^{T}
\mathbf{X}_{2},\mathbf{v}^{T}\mathbf{X}_{3}\leq
\mathbf{v}%
^{T}\mathbf{X}_{2}\bigr)-\frac{1}{4\lambda(1-\lambda)},\nonumber
\end{eqnarray}
where $\mathbf{X}_{1},\mathbf{X}_{2},\mathbf{X}_{3}$ are
i.i.d. from
the common DF.
\end{pf*}
\begin{pf*}{Proof of Theorem~\ref{Thm-Consistent&Monotone}}
Suppose that $\mathbf{X}\prec_{\lst }\mathbf{Y}$. Then for some
$\mathbf{s}_{\ast}\in\mathcal{S}_{+}^{p-1}$ we have $\mathbf
{s}_{\ast
}^{T}\mathbf{X}\prec_{\st }\mathbf{s}_{\ast}^{T}\mathbf{Y} $ which
implies that $\mathbb{P} ( \mathbf{s}_{\ast}^{T}\mathbf{X}%
\leq\mathbf{s}_{\ast}^{T}\mathbf{Y} ) >1/2$. By definition
$\mathbb{P} ( \mathbf{s}_{\max}^{T}\mathbf{X}\leq\mathbf
{s}%
_{\max}^{T}\mathbf{Y} ) \geq\mathbb{P} ( \mathbf
{s}_{\ast
}^{T}\mathbf{X}\leq\mathbf{s}_{\ast}^{T}\mathbf{Y} ) $ so
$\Psi(\mathbf{s}_{\max})>1/2$. It follows from the proof of Theorem
\ref{Them-LST1} that $\Psi_{n,m}(\widehat{\mathbf
{s}%
}_{\max})\rightarrow\Psi(\mathbf{s}_{\max})$ with probability one.
Thus,%
\[
S_{n,m}=N^{1/2} \bigl( \Psi_{n,m}(\widehat{
\mathbf{s}}_{\max
})-1/2 \bigr) \stackrel{\mathrm{a.s.}} {\longrightarrow}\infty
\qquad\mbox{as }n,m\rightarrow\infty.
\]
Therefore by Slutzky's theorem,%
\[
\mathbb{P} ( S_{n,m}>q_{n,m,1-\alpha};H_{1} ) \rightarrow1
\qquad\mbox{as }n,m\rightarrow\infty,
\]
where $q_{n,m,1-\alpha}$ is the critical value for an $\alpha$ level test
based on samples of size $n$ and $m$ and $q_{n,m,1-\alpha}\rightarrow
q_{1-\alpha}$. Hence the test based on $S_{n,m}$ is consistent. Consistency
for $I_{n,m}$ and $I_{n,m}^{+}$ is established in a similar manner.

Now assume that $\mathbf{X}\preceq_{\lst }\mathbf{Y}\preceq_{\lst }\mathbf
{Z}$ so that $\mathbf{s}^{T}\mathbf{Y}\preceq_{\st }\mathbf{s}^{T}\mathbf{Z}$
for all $\mathbf{s}\in
\mathcal{S}_{+}^{p-1}$. Fix $\mathbf{x}_{i}$, $i=1,\ldots,n$, and choose
$\mathbf{s}\in\mathcal{S}_{+}^{p-1}$. Without any loss of generality
assume that $\mathbf{s}^{T}\mathbf{x}_{1}\leq$ $\mathbf{s}%
^{T}\mathbf{x}_{2}\leq\cdots\leq\mathbf{s}^{T}\mathbf{x}_{n}$.
Define $U_{j}=\sum_{i=1}^{n}\mathbb{I}_{(\mathbf{s}^{T}\mathbf
{x}%
_{i}\leq\mathbf{s}^{T}\mathbf{Y}_{j})}$ and $V_{j}=\sum_{i=1}%
^{n}\mathbb{I}_{(\mathbf{s}^{T}\mathbf{x}_{i}\leq\mathbf{s}%
^{T}\mathbf{Z}_{j})}$. Clearly $U_{j}$ and $V_{j}$ take values in
$J= \{ 0,\ldots,n \}$. Now, for $k\in J$ we have
\[
\mathbb{P} ( U_{j}\geq k ) =\mathbb{P} \bigl( \mathbf{s}%
^{T}
\mathbf{Y}_{j}\geq\mathbf{s}^{T}\mathbf{x}_{k} \bigr)
\leq\mathbb{P} \bigl( \mathbf{s}^{T}\mathbf{Z}_{j}\geq\mathbf
{s}%
^{T}\mathbf{x}_{k} \bigr) =\mathbb{P} (
V_{j}\geq k),
\]
where we use the fact that $\mathbf{s}^{T}\mathbf{Y}\preceq_{\st }\mathbf
{s}^{T}\mathbf{Z.}$ It follows that $U_{j}\preceq_{\st }V_{j}$ for
$j=1,\ldots,m$. Moreover $ \{ U_{j} \} $ and
$ \{ V_{j} \} $ are all independent and it follows from Theorem
1.A.3 in \citet{ShaSha07} that $\sum_{j=1}^{m}U_{j}\preceq
_{\st }\sum_{j=1}^{m}V_{j}$. Thus $\sum_{i=1}^{n}\sum_{j=1}^{m}\mathbb{I}%
_{(\mathbf{s}^{T}\mathbf{x}_{i}\leq\mathbf{s}^{T}\mathbf
{Y}%
_{j})}\preceq_{\st }\sum_{i=1}^{n}\sum_{j=1}^{m}\mathbb{I}_{(\mathbf
{s}%
^{T}\mathbf{x}_{i}\leq\mathbf{s}^{T}\mathbf{Z}_{j})}$. The
latter
holds for every value of $\mathbf{x}_{1},\ldots,\mathbf
{x}_{n}$, and
therefore it holds unconditionally as well, that is,%
\[
\sum_{i=1}^{n}\sum
_{j=1}^{m}\mathbb{I}_{(\mathbf{s}^{T}\mathbf
{X}%
_{i}\leq\mathbf{s}^{T}\mathbf{Y}_{j})}
\preceq_{\st }\sum_{i=1}^{n}
\sum_{j=1}^{m}
\mathbb{I}_{(\mathbf{s}^{T}\mathbf{X}_{i}\leq
\mathbf{s}^{T}\mathbf{Z}_{j})}.
\]
It follows that $\Psi_{n,m}^{\mathbf{X},\mathbf{Y}}(\mathbf
{s}%
)\preceq_{\st }\Psi_{n,m}^{\mathbf{X},\mathbf{Z}}(\mathbf{s})$ for all
$\mathbf{s}\in\mathcal{S}_{+}^{p-1}$ where $\Psi_{n,m}^{\mathbf{X},\mathbf
{Y}}(\mathbf{s})$ and
$\Psi_{n,m}^{\mathbf{X},\mathbf{Z}}(\mathbf{s})$ are defined in
(\ref{Psi-nm}) and the superscripts emphasize the different arguments
used to evaluate them.
Thus%
\[
\Psi_{n,m}^{\mathbf{X},\mathbf{Y}}\bigl(\widehat{\mathbf{s}}^{\mathbf{X},\mathbf{Y}}_{\max
} \bigr)\preceq_{\st }\Psi_{n,m}^{\mathbf{X},\mathbf {Z}}\bigl(\widehat{
\mathbf{s}}^{\mathbf{X},\mathbf{Y}}_{\max}\bigr)\preceq_{\st }
\Psi_{n,m}^{\mathbf{X},\mathbf{Z}}\bigl(\widehat{\mathbf{s}}^{\mathbf{X},\mathbf{Z}}_{\max}
\bigr)
\]
and as a consequence $\mathbb{P}(S_{n,m}^{X,Y}>q_{n,m,1-\alpha})\leq
\mathbb{P}(S_{n,m}^{X,Z}>q_{n,m,1-\alpha})$ as required. The
monotonicity of
the power function of $I_{n,m}$ and $I_{n,m}^{+}$ follows immediately
from the
fact that $\Psi_{n,m}^{\mathbf{X},\mathbf{Y}}(\mathbf{s}%
)\preceq_{\st }\Psi_{n,m}^{\mathbf{X},\mathbf{Z}}(\mathbf
{s})$ for
all $\mathbf{s}\in\mathcal{S}_{+}^{p-1}$.
\end{pf*}
\end{appendix}

\section*{Acknowledgments}

We thank Grace Kissling (NIEHS), Alexander Goldenshluger, Yair Goldberg
and Danny Segev (University of Haifa), for their useful comments and
suggestions. We also thank the Editor, Associate Editor and two referees for their
input which improved the paper.



\printaddresses

\end{document}